\documentclass[12pt, reqno]{amsart}
\usepackage{amsmath, amsthm, amscd, amsfonts, amssymb, graphicx, color}
\usepackage[bookmarksnumbered, colorlinks, plainpages]{hyperref}
\hypersetup{colorlinks=true,linkcolor=red, anchorcolor=green, citecolor=cyan, urlcolor=red, filecolor=magenta, pdftoolbar=true}

\newtheorem{theorem}{Theorem}[section]
\newtheorem{lemma}[theorem]{Lemma}
\newtheorem{proposition}[theorem]{Proposition}
\newtheorem{question}[theorem]{Question}
\newtheorem{corollary}[theorem]{Corollary}
\theoremstyle{definition}

\newtheorem{example}[theorem]{Example}

\theoremstyle{remark}
\newtheorem{remark}[theorem]{Remark}
\numberwithin{equation}{section}

\begin{document}

\setcounter{page}{1}

\title[]{A correction of ``D. Taghizadeh, M. Zahraei, A. Peperko and N. H. Aboutalebi. \textrm{On the numerical ranges of matrices in max algebra,
} Banach J. Math. Anal. \textbf{14} (2020),  1773--1792.''}

\author[Thaghizadeh, Zahraei, Haj Aboutalebi, Peperko,  Fallat, Khorami]{Davod Thaghizadeh$^{1},$  Mohsen Zahraei$^{1*},$ Narges Haj Aboutalebi$^{2},$ Aljo\v{s}a Peperko$^{3, 4},$    Shaun Fallat $^{5}$ \MakeLowercase{and} Reza Tayebi Khorami$^{1}$ }

\address{$^{1}$Department of Mathematics, Ahvaz Branch, Islamic Azad University, Ahvaz, Iran.}
\email{\textcolor[rgb]{0.00,0.00,0.84}{t.davood1411@gmail.com
}}

\address{
$^{1*}$
Department of Mathematics, Ahvaz Branch, Islamic Azad University, Ahvaz, Iran.}
\email{\textcolor[rgb]{0.00,0.00,0.84}{mzahraei326@gmail.com 
}}

\address{$^{2}$Department of Mathematics, Shahrood Branch, Islamic Azad University, Shahrood, Iran.
}
\email{\textcolor[rgb]{0.00,0.00,0.84}{aboutalebi.n@yahoo.com
}}

\address{$^{3}$ Faculty of Mechanical Engineering, University of Ljubljana, Aškečeva 6, SI-1000 Ljubljana, Slovenia. }
\email{\textcolor[rgb]{0.00,0.00,0.84}{aljosa.peperko@fs.uni-lj.si
}}

\address{
$^{4}$Institute of Mathematics, Physics and Mechanics Jadranska 19, SI-1000 Ljubljana, Slovenia. }
\email{\textcolor[rgb]{0.00,0.00,0.84}{aljosa.peperko@fmf.uni-lj.si
}}

\address{$^{5}$ Department of Mathematics and Statistics, University of Regina, Saskatchewan, Canada.}
\email{\textcolor[rgb]{0.00,0.00,0.84}{shaun.fallat@uregina.ca
}}

\address{$^{1}$Department of Mathematics, Ahvaz Branch, Islamic Azad University, Ahvaz, Iran.}
\email{\textcolor[rgb]{0.00,0.00,0.84}{
r.t.khorami@gmail.com  }}






\date{Received: xxxxxx; Revised: yyyyyy; Accepted: zzzzzz.
\newline \indent $^{*}$Corresponding author}

\begin{abstract}
We correct some unfortunate mistakes that appeared in the article D. Taghizadeh, M. Zahraei, A. Peperko and N. H. Aboutalebi, \textit{On the numerical ranges of matrices in max algebra}, Banach J. Math. Anal.,  \textbf{14} (2020), pp. 1773--1792 concerning certain notions of the numerical range in the max algebra setting. To do this we also include a study of the characteristic max polynomial and correspondingly the max $k$-spectrum and the $k$-tropical spectrum. We also pose a nonresolved open question.


\end{abstract}
 \maketitle
 
 \section{{\bf Introduction}}
 
 In \cite{TZAH} different notions of numerical ranges in max algebra were studied. Unfortunately, some errors were identified in the original proofs of a few results (for example, in Sections 4 and 5 of \cite{TZAH}). The aim of this work is to correct and rectify these lamentable mistakes and provide more clarity on both their verification and corresponding applications.
 
Suppose that $A=(a_{ij})\in M_n(\mathbb{R}_+)$ and let $1\leq k\leq n$ be a positive integer.  It was stated in \cite[Theorem 3]{TZAH} that the $k-$numerical range $W_{\max}^k(A)$ in max algebra equals
 \begin{equation}\label{TO_CORRECT}
 W_{\max}^k(A)=[c,d],
 \end{equation}
 where $c=\min\{\displaystyle\oplus_{j=1}^ka_{i_ji_j}:~1\leq i_1<i_2<\cdots <i_k\leq n\}$ and $d=\displaystyle\max_{1\leq i,j\leq n}a_{ij}$. Although as written the statement \cite[Theorem 3]{TZAH} is correct for $1 \le k < n$ (see \cite[Theorem 1]{TZAH} and Theorem \ref{main_correction} below), it turns out that it is not correct for $k=n$. This inaccuracy was observed by S. Gaubert  and was communicated to the authors in an editorial communication \cite{SG}. S. Gaubert observed that in fact
 \begin{equation*}
W_{max}^{n}(A)=\{\displaystyle\max_{1\leq i\leq n}a_{ii} \},
\end{equation*}
contradicting the statement that the upper bound of $W_{max}^{k}(A)$ is $d =  \displaystyle\max_{1\leq i, j \leq n}a_{ij},$ as claimed in \cite[Theorem 3]{TZAH}. The authors of the current article apologize to the readers for this mistake and also for some other mistakes from \cite{TZAH}, which we correct in this article. To ensure that the results presented are clear and precise, we first list which results from \cite{TZAH} are correct (including their proofs) and which need to be revised and clarified.

The following parts of the article \cite{TZAH} are correct (including their proofs) as listed:
\begin{itemize}
\item \cite[Sections 2 and 3]{TZAH} ;

\item \cite[Equations (14) and (15), Remark 5, Propositions  10 and 11, Example 4 from Section 4]{TZAH};

\item \cite[Remark 7, Theorem 5(vii) from Section 5]{TZAH}.
 \end{itemize}
 The following parts of the article \cite{TZAH} contain mistakes:
\begin{itemize}
\item \cite[Remark 4, Theorem 3 for $k=n$, Example 3 and Proposition 7 for $k=n$ in Section 4]{TZAH}, the corrections are given in [Remark \ref{Remark3.5LocallyLip}, Theorem \ref{Main Theorem}, Remark \ref{Exaple3TZAH} and Remark \ref{proposition7TZH}]  below.

\item \cite[discussion after Definition 6 and before Remark 7 in Section 5]{TZAH}, the corrections are given in [discussion after Equality (\ref{druga})] below.
\end{itemize}
The following parts of the article \cite{TZAH} are correct, but require different or ammended proofs (which we provide in the current article):
\begin{itemize}
\item \cite[Theorem 3 in the case $1 \le k <n$]{TZAH} (see Theorem \ref{main_correction} below);
\item \cite[Theorem 4]{TZAH} (see Theorem \ref{TheoremFirstknumer} below); 
\item \cite[Proposition 7 in the case $1 \le k <n$]{TZAH} (see Proposition \ref{inclusionsigmak} and Theorem \ref{trop_cor} below);
\item \cite[Propositions 8 and 9]{TZAH} (they follow from Proposition \ref{inclusionsigmak} and Theorems \ref{Main Theorem} and \ref{trop_cor} below);
\item \cite[Example 5 from Section 5]{TZAH} (see Remark \ref{extend_C numerical range} below);
\item \cite[Theorem 5(i)-(iv), (vi) and Corollary 4 from Section 5]{TZAH} (see Theorem \ref{THCNumerical} below).
\end{itemize}
It is not clear if the set  $\mathcal{X}_{n \times k} $ as defined in \cite[Remark 4]{TZAH} is a connected set for $1 < k < n$ as it was stated  in \cite[Remark 4]{TZAH}. This is left as an open question (see Question \ref{open_qu1}).

 
  \section{{\bf Preliminaries}}
  
  A conventional max  algebra consists of the set of nonnegative real numbers equipped with the basic operations of multiplication  $a\otimes b=ab,$ and maximization  $a\oplus b=\max\{a,b\}$ (see also e.g. \cite{M.PEPERKO1}, \cite{GMW}, \cite{Butkovic}, \cite{BCOQ}, \cite{KM} \cite{TS} and the references cited within).   
  For  $A=(a_{ij}) \in M_{m\times n}(\mathbb{R}),$  we say that $A$  is positive (nonnegative) and write $A>0$ $(A \geq 0)$  if $a_{ij}>0$ $(a_{ij} \geq 0)$  for  $1 \leq i \le m$, $1\le j \leq n.$  Let $\mathbb{R}_{+}$ be the set of all nonnegative real numbers and $M_{m \times n}(\mathbb{R}_{+})$ denote the set of all $m \times n$ nonnegative (real) matrices.
The notions $M_{n}(\mathbb{R}_{+})$ and $\mathbb{R}_{+}^{n}$ are abbreviations for $M_{n \times n}(\mathbb{R}_{+})$  and  $M_{n \times 1}(\mathbb{R}_{+}),$  respectively.  

Let  $A=(a_{ij})\in M_{m \times n}(\mathbb{R_{+}})$  and  $B=(b_{ij}) \in M_{n \times l}(\mathbb{R}_{+}).$  The  product of $A$  and  $B$  in  a max algebra is denoted  by  $A\otimes B,$   and for  $1 \leq i \le m$, $1\le j \leq l$, $(A\otimes B)_{ij}=\displaystyle\max_{k=1,\ldots,n} a_{ik}b_{kj}.$ In particular,  for $x \in \mathbb{R}_{+}^{n}$ the vector $A \otimes x$ is defined by $(A\otimes x)_i = \displaystyle\max _{k=1,\ldots,n} a_{ik}x_k $ for $i=1, \ldots , m $. 
  If $A, B \in M_{n}(\mathbb{R}_{+})$, then the max sum $A\oplus B$ in a max algebra is defined by $(A\oplus B)_{ij}=\displaystyle\max \{a_{ij}, b_{ij}\}$ for $i,j=1, \ldots , n$.
  The notation  $A_{\otimes}^{2}$  refers to  $A\otimes A,$  and  $A_{\otimes}^{k}$  denotes  the $k$th  power  of $A$ in a max algebra.
For $A \in M_{ n}(\mathbb{R})$ and $x \in \mathbb{R}^{n},$ let $\displaystyle \|A\|= \max _{i,j=1, \ldots , n} |a_{ij}|$, $tr_{\otimes}(A)=\max _{i=1, \ldots , n} a_{ii}$ and $\|x\|=\displaystyle\max _{i=1, \ldots , n} |x_i|$. Finally, let $A^t$ and $x^t$ denote the transpose of matrix $A$, and the vector $x$, respectively.

In \cite{TS}, the numerical range of a given square matrix was introduced and described in the setting of the max-plus algebra, while its isomorphic version in a max algebra setting was studied in \cite{TZAH}. 
If  $A\in M_n(\mathbb{R_+})$, then, the max numerical range $W_{max}(A)$ of $A$, is defined by 
$${W_{\max }}(A) = \left\{ {{x^t} \otimes A \otimes x:\,\,\,x \in {\mathbb{R} }_{+}^n,\,{x^t} \otimes x = 1\,} \right\}.$$

The following theorem  was proved in \cite[Theorem 3.7]{TS} an alternative proof was given in \cite[Theorem 2 and Remark 2]{TZAH}.

\begin{theorem}\label{wmaxA}
Let  $A=(a_{ij})\in M_n(\mathbb{R_+})$ be a nonnegative matrix. Then
\[
W_{\max }(A)=[a, b] \subseteq \mathbb{R_{+}}, 
\]
where $a=\displaystyle\min_{1 \leq i \leq n}a_{ii}$ and $b=\displaystyle \max_{1 \leq i,j \leq n}a_{ij}=\|A\|.$    
\end{theorem}
 Let $U\in M_n(\mathbb{R}_{+})$ and let  $I_{n}$ denote  the $n \times  n $
identity matrix. If $U^t\otimes U=U\otimes U^t=I_n,$ then $U$ is called {\em unitary} in   a max algebra and we  denote the collection of all unitary matrices by
$$\mathcal{U}_{n}=\{U \in M_n(\mathbb{R}_{+}): U^{t}\otimes U=U\otimes U^t= I_n \}.$$
It is well known (see e.g. \cite{Butkovic}) that 
$A\in M_n(\mathbb{R}_+)$ is unitary in   a max algebra if and only if $A$ is  a permutation matrix. Thus the set $\mathcal{U}_{n}$ is the group of permutation matrices.

 Let $k$ and $n$ be  positive integers such that $k\leq n$ and $A\in M_n(\mathbb{R}_+)$.   
 A set $\mathcal{X}_{n \times k} \subset 
M_{n \times k}(\mathbb{R}_{+})$ is defined by
$$\mathcal{X}_{n \times k}= \{X \in M_{n \times k}(\mathbb{R}_{+}): X^{t}\otimes X = I_{k}\}.$$
It is known that 
 for the case $k = n,$  $\mathcal{X}_{n \times n}$ coincides with $\mathcal{U}_{n}$ (see, e.g., \cite[Lemma 4.84]{BCOQ} or \cite[Corollary 2.4]{KM}, and also see Lemma \ref{properties} below). 

For $X \in M_{n \times k}(\mathbb{R}_{+})$ let $x^{(i)}$ denote the $i$th column of $X$ for $i=1,\ldots , k$. So for $X=(x_{ji})_{j=1, \ldots, n} ^{i=1, \ldots ,k}$ we have $x^{(i)}_j =x_{ji}$. Let $A \in M_{n}(\mathbb{R}_{+})$ and   $1 \le k \leq n$.   The {\em max $k$-numerical range $W_{\max}^k(A)$ of $A$} in a max algebra was introduced in \cite[Section 4]{TZAH} and is defined by
%
\begin{eqnarray*}
W_{max}^{k}(A) &=&\{\bigoplus_{i=1}^{k} \displaystyle (x^{(i)})^{t}\otimes A \otimes x^{(i)}:  ~X=[x^{(1)}, x^{(2)}, \ldots, x^{(k)}] \in \mathcal{X}_{n \times k} \}\\
&=& \{tr_{\otimes}(X^{t}\otimes A \otimes X):~ X=[x^{(1)}, x^{(2)}, \ldots, x^{(k)}] \in \mathcal{X}_{n \times k}   \}. 
\end{eqnarray*}
Note that $W_{max}^{1}(A)=W_{max}(A)$ 
and 
\[
tr_{\otimes}(X^{t}\otimes A\otimes X)=\displaystyle (x^{(1)})^{t}\otimes A\otimes x^{(1)}\oplus (x^{(2)})^{t}\otimes A\otimes x^{(2)}\oplus \cdots \oplus (x^{(k)})^{t}\otimes A\otimes x^{(k)},
\]
for all $ X=[x^{(1)}, x^{(2)}, \ldots, x^{(k)}] \in M_{n \times k}(\mathbb{R}_{+})$. Observe  that for any $X \in  \mathcal{X}_{n \times k}$ it follows that
 \begin{equation}
   (x^{(i)})^{t}\otimes x^{(j)}=\delta_{ij}= 
  \begin{cases}
  1  &    i=j, \\
  0  &    i \neq j.
  \end{cases}
  \label{10}
  \end{equation}


\smallskip

\section{\textbf{Corrections from \cite[Section 4]{TZAH}}}

It was stated in \cite[Remark 4 and Theorem 3]{TZAH} that for $A \in M_{n}(\mathbb{R}_{+})$ and $1 \le k \le n$ the set $\mathcal{X}_{n \times k}$ is a connected set and that
\begin{equation} 
  W_{\max }^{k}(A)=[c, d],
\label{c,d}
  \end{equation}
  where $\displaystyle c=\min \{\bigoplus_{j=1}^{k} a_{i_{j}i_{j}}:  1 \leq i_{1}< i_{2} < \cdots < i_{k} \leq n  \}$  and  
  $\displaystyle d=\max_{1 \leq i, j\leq n}a_{ij}.$
  However, as justly pointed out by S. Gaubert in an editorial communication \cite{SG} for $k=n$ the above statement is not true in general.
  The set $\mathcal{X}_{n \times n}=\mathcal{U}_n$ is not connected (being a group of permutation matrices) and the equality
   \begin{equation}\label{WmaxNnumberone}
W_{max}^{n}(A)=\{\displaystyle\max_{1\leq i\leq n}a_{ii} \}
\end{equation}
  holds.
  
   To establish (\ref{WmaxNnumberone}) observe that
  \[
\displaystyle\max_{1\leq i \leq n}(x^{(i)})^{t}\otimes A \otimes x^{(i)}=\max_{1\leq i \leq n}a_{\sigma(i), \sigma(i)},
\]
where $\sigma$ is the permutation  represented by the matrix $X\in \mathcal{U}_{n}$. Hence 
\begin{equation*}
W_{max}^{n}(A)=\{ \displaystyle\max_{1\leq i\leq n}a_{\sigma(i), \sigma(i)}: \sigma\in \sigma_{n}\}=\{\displaystyle\max_{1\leq i\leq n}a_{ii} \},
\end{equation*}
which proves (\ref{WmaxNnumberone}).
Note that $\sigma_{n}$  denotes the symmetric group (group of permutations) on $\{1,2,\ldots,n\}$. 

\begin{remark}
The permutation group $\mathcal{U}_{n}$ has exactly $n!$ elements and for each two different permutation matrices $V_{1}, V_{2}\in \mathcal{U}_{n},$ 
$
\Vert V_{1}-V_{2} \Vert=1.
$
So it follows that
\[
\mathcal{U}_{n}=\displaystyle  \bigcup_{i=1}^{n!}B(V_{i}; \frac{1}{2}),~ V_{i}\in \mathcal{U}_{n}~ i=1, 2, \ldots, n!,
\] 
where $\displaystyle\{B(V_{i}; \frac{1}{2}) \}_{i=1}^{n!}$ is a collection of  disjoint  non empty open balls of  $\mathcal{U}_{n}$.
The well known above argument verifies that $\mathcal{U}_{n}$ is totally disconnected. Moreover, $\mathcal{U}_{n}$ is compact and Hausdorff.
\end{remark}

As seen above the equality in (\ref{c,d}) need not hold for $k=n$. On the other hand, for the case when $k=1$  (\ref{c,d}) holds by Theorem \ref{wmaxA}, and clearly the set $ \mathcal{X}_{n \times 1}$ is a compact connected set. We establish below, in Theorem  \ref{main_correction}, that (\ref{c,d}) holds whenever $1<k<n$. Moreover, our proof of Theorem  \ref{main_correction} is constructive in nature.

It remains however unclear  if the set $ \mathcal{X}_{n \times k}$ is  connected when $1 < k < n$, which leads to the query below.

\begin{question} Is the set $\mathcal{X}_{n \times k} $ a connected set for $1 < k < n$?
\label{open_qu1}
\end{question}

Next we consider the following lemma regarding the set $\mathcal{X}_{n \times k}$.
\begin{lemma}
Let $1 \le k \le n$ and $X \in \mathcal{X}_{n \times k}$. Then the following properties hold.

(i) For each $i \neq j$, $i,j \in \{1, \ldots , k \}$ and each $l \in \{1, \ldots , n\}$ either $x^{(i)}_{l}=0$ or $x^{(j)}_{l}=0$.

(ii) For each $l \in \{1, \ldots , n\}$ and $i \in \{1, \ldots , k\}$, $x^{(i)}_{l} \in [0,1]$.

(iii) For each  $i \in \{1, \ldots , k\}$ there exists $l \in \{1, \ldots , n\}$ such that $x^{(i)}_{l}=1$ and that $x^{(j)}_{l}=0$ for all $j\neq i$, $j \in \{1, \ldots , k\}$.

(iv) $X$ has a $k \times k$ permutation submatrix (that is, $X$ has a $k \times k$ submatrix that is a permutation matrix).

(v) $\mathcal{X}_{n \times n}= \mathcal{U}_n$.
\label{properties}
\end{lemma}
\begin{proof} Property (i) follows from (\ref{10}). Properties (ii) and (iii) follow from (i) and (\ref{10}). Property (iv) follows from (i) and (iii), while property (v) follows from (iv).
\end{proof}
%
Most of the  following result was stated in \cite[Theorem 4]{TZAH}. Since the original proof contained errors, we provide a reformed and valid proof below.
 \begin{theorem}\label{TheoremFirstknumer}
  Let $A\in M_{n}(\mathbb{R}_{+})$  and let $1\leq k \leq n$ be a positive integer. Then the following assertions  hold.
 \begin{itemize}
 \item[(i)]  $W_{max}^{k}(\alpha A\oplus \beta I)=\alpha W_{max}^{k}(A)\oplus \beta,$  $W_{max}^{k}(A\oplus B)\subseteq W_{max}^{k}(A)\oplus W_{max}^{k}(B)$  and $W_{max}^{n}(A\oplus B)=W_{max}^{n}(A)\oplus W_{max}^{n}(B),$  where $\alpha, \beta\in \mathbb{R}_{+}$  and $B\in M_{n}(\mathbb{R}_{+})$;
 
 \item[(ii)] $W_{max}^{k}(U^{t}\otimes A\otimes U)=W_{max}^{k}(A)$ if $U\in \mathcal{U}_{n}$;
 
 \item[(iii)]  If $B\in M_{m}(\mathbb{R}_{+})$ is a principal submatrix of $A$   and $k\leq m,$  then $W_{max}^{k}(B)\subseteq W_{max}^{k}(A)$. 
 Consequently, if $V=[e_{i_{1}}, e_{i_{2}}, \ldots, e_{i_{s}}]\in M_{n\times s}(\mathbb{R}_{+}),$   where $1\leq k\leq s \leq n,$  then
 $W_{max}^{k}(V^{t}\otimes A\otimes V)\subseteq W_{max}^{k}(A),$  and equality holds if $s=n$, and where $e_l$ denotes the standard basis vector in $\mathbb{R}^n$;
 
 \item[(iv)] $W_{max}^{k}(A^{t})=W_{max}^{k}(A)$;
 
 \item[(v)] If $k<n,$  then $W_{max}^{k+1}(A)\subseteq W_{max}^{k}(A)$. Consequently, 
\[
 W_{max}^{n}(A)\subseteq W_{max}^{n-1}(A)\subseteq \cdots \subseteq W_{max}^{2}(A)\subseteq W_{max}(A).
 \]
 \end{itemize}
 \end{theorem}
 \medskip
 
 \begin{proof}
 \begin{itemize}
  \item[(i)] Let $z\in W_{max}^{k}(\alpha A\oplus \beta I)$. So $z=\displaystyle\max_{1\leq i \leq k}(x^{(i)})^{t}\otimes (\alpha A\oplus \beta I)\otimes x^{(i)}$  for some $X\in \mathcal{X}_{n\times k}$ and hence $z=\alpha \displaystyle(\max_{1\leq i \leq k}(x^{(i)})^{t}\otimes A \otimes x^{(i)})\oplus \beta$. This implies that $z\in \alpha W_{max}^{k}(A)\oplus \beta$. For the reverse inclusion, let $z\in \alpha W_{max}^{k}(A)\oplus \beta$. So $z=\alpha \displaystyle(\max_{1\leq i \leq k}(x^{(i)})^{t}\otimes A \otimes x^{(i)})\oplus \beta$ for some $X\in \mathcal{X}_{n\times k}$ and it follows that $z\in W_{max}^{k}(\alpha A\oplus \beta I)$.

  For the second part, suppose that $z\in W_{max}^{k}(A\oplus B)$. Then \begin{eqnarray*}
  z &=&\displaystyle \max_{1\leq i\leq k}(x^{(i)})^{t}\otimes(A\oplus B)\otimes x^{(i)}\\
  &=&\displaystyle \max_{1\leq i\leq k}(((x^{(i)})^{t}\otimes A\otimes x^{(i)})\oplus ((x^{(i)})^{t}\otimes B\otimes x^{(i)}))\\
  &=& \displaystyle \max_{1\leq i\leq k}((x^{(i)})^{t}\otimes A\otimes x^{(i)})\oplus \displaystyle \max_{1\leq i\leq k} ((x^{(i)})^{t}\otimes B\otimes x^{(i)}).
 \end{eqnarray*} 
  This shows that $z\in W_{max}^{k}(A)\oplus W_{max}^{k}(B)$. The equality in the case $k=n$ follows from (\ref{WmaxNnumberone}).
  \item[(ii)]
  Let $z\in W_{max}^{k}(U^{t}\otimes A\otimes U)$  where $U\in \mathcal{U}_{n}$. Then
  \[
  z=\displaystyle\max_{1\leq i \leq k}(U\otimes x^{(i)})^{t}\otimes A\otimes (U\otimes x^{(i)})~~\mathrm{for~ some}~ X=[x^{(1)}, x^{(2)}, \ldots, x^{(k)}]\]
 is in $\mathcal{X}_{n\times k}.$
  Since $U\otimes X=[U\otimes x^{(1)}, U\otimes x^{(2)}, \ldots, U\otimes x^{(k)}]\in \mathcal{X}_{n\times k}$  it follows that $z\in W_{max}^{k}(A)$.

To prove the reverse inclusion observe that for $X\in  \mathcal{X}_{n\times k}$ and $U\in \mathcal{U}_{n}$ it holds that $U^t \otimes X \in \mathcal{X}_{n\times k}$
and 
$X^t \otimes A \otimes X= (U^t \otimes X)^t \otimes U^t \otimes A \otimes U \otimes  (U^t \otimes X), $
and so $W_{max}^{k}(A) \subset W_{max}^{k}(U^{t}\otimes A\otimes U)$.

  
  \item[(iii)] By (ii) we may assume, without loss of generality,  that
  \begin{equation*}\label{B_Corner}
  A=\left [
\begin{array}{cc}
B  & \star\\
\star & \star
\end{array}
\right].
  \end{equation*}
  Let $z\in W_{max}^{k}(B)$. So $z=\displaystyle\max_{1\leq i \leq k} (x^{(i)})^{t}\otimes B\otimes x^{(i)}$  for some $X=[x^{(1)}, x^{(2)}, \ldots, x^{(k)}]\in \mathcal{X}_{m\times k}$. Now by taking 
  
\begin{math}
Y=\left[
\begin{array}{cccc}
x^{(1)} &  x^{(2)} & \ldots &  x^{(k)}\\
 O_{(n-m)\times 1} &  O_{(n-m)\times 1} & \ldots & O_{(n-m)\times 1}
\end{array}
\right]\in \mathcal{X}_{n\times k}
\end{math}
($O_{(s)\times 1}$ is the zero $s$-vector)
we have 
\[
z=\displaystyle\max_{1\leq i\leq k}(y^{(i)})^{t}\otimes A\otimes y^{(i)}\in W_{max}^{k}(A).
\]
For the second part, suppose that $z\in W_{max}^{k}(V^{t}\otimes A\otimes V)$, where $V=[e_{i_{1}}, e_{i_{2}}, \ldots, e_{i_{s}}]$. Then 
\[
z=\displaystyle\max_{1\leq i \leq k}(V\otimes x^{(i)})^{t}\otimes A\otimes (V\otimes x^{(i)})~~\mathrm{for~ some}~ X=[x^{(1)}, x^{(2)}, \ldots, x^{(k)}] \]
is in $\mathcal{X}_{s\times k}$.
Since $V\otimes X=[V\otimes x^{(1)}, V\otimes x^{(2)}, \ldots, V\otimes x^{(k)}]\in \mathcal{X}_{n\times k},$ $z\in W_{max}^{k}(A)$ and so $W_{max}^{k}(V^{t}\otimes A\otimes V)\subseteq W_{max}^{k}(A)$. For the case when $s=n$ equality follows from 
(ii).

\item[(iv)]
Since for each $X \in  \mathcal{X}_{n\times k}$ we have 
$tr_{\otimes}(X^{t}\otimes A \otimes X)=tr_{\otimes}(X^{t}\otimes A^t \otimes X), $
it follows that
\[
W_{max}^{k}(A^{t})=W_{max}^{k}(A).
\]
\item[(v)] Let $z \in W_{\max }^{k+1}(A)$. So there exist  $X=[x^{(1)}, x^{(2)}, \ldots, x^{(k)}, x^{({k+1})}] \in \mathcal{X}_{n \times (k+1)}$ such that 
\[
z=\bigoplus_{i=1}^{k+1} (x^{(i)})^{t}\otimes A \otimes x^{(i)}.
\]
Now by (ii),   assume, without loss of generality,
\[
 (x^{(1)})^{t}\otimes A \otimes x^{(1)}=\displaystyle\min_{1\leq i\leq k+1}(x^{(i)})^{t}\otimes A\otimes x^{(i)}.
 \]
 Hence, by setting $Y=[ x^{(2)}, \ldots, x^{(k)}, x^{(k+1)}], $  we have $Y \in \mathcal{X}_{n \times k}$ and hence
 \[
 z=\bigoplus_{i=2}^{k+1} (x^{(i)})^{t}\otimes A \otimes x^{(i)}.
 \]
 This implies that $z \in W_{\max }^{k}(A) $, which completes the proof.
 \end{itemize}
\end{proof} 

It was noted in \cite[Remark 4]{TZAH} that for $1\leq k\le n$ the set $\mathcal{X}_{n\times k}$ is compact and that the mapping $f_{A}: \mathcal{X}_{n\times k}\longrightarrow \mathbb{R}_{+}$ defined by
 $$f_{A}(X):=tr_{\otimes}(X^t\otimes A\otimes X)$$
 is locally Lypschitz. These claims are correctly verified below.
\begin{remark}\label{Remark3.5LocallyLip} 
Let $k$ and $n$ be  positive integers such that $1\leq k \le n$ and $A\in M_n(\mathbb{R}_+)$. Suppose $\{X_{m}\}$ is a given sequence in $\mathcal{X}_{n\times k}$ such that $X_{m}\longrightarrow X$ as $m\longrightarrow \infty$.     By  continuity of  the max product we have
\[
X_{m}^{t}\otimes X_{m}\longrightarrow X^{t}\otimes X
\]
 and since $X_{m}^{t}\otimes X_{m}=I_{k}$ for all $m=1, 2, \ldots,$ it follows that $X^{t}\otimes X=I_{k}$. This verifies that 
 $X\in \mathcal{X}_{n\times k}$ and thus $\mathcal{X}_{n\times k}$ is closed. Since $\mathcal{X}_{n\times k}$ is also a bounded set,   $\mathcal{X}_{n\times k}$ is a compact set.
\end{remark}

For any $n \times n$ real matrix $A=(a_{ij})$,  the absolute value of $A$ is defined entrywise $|A| = (|a_{ij}|)\in M_{n}(\mathbb{R}_+) $.

\begin{lemma}\label{lemma1}
Let $A=(a_{ij}), B=(b_{ij})\in M_n(\mathbb{R}_+)$. Then
\begin{equation}
\label{Lip_of_tr}
\vert tr_{\otimes}(A) - tr_{\otimes}(B) \vert\leq tr_{\otimes}\vert A- B \vert \le
 \|A-B\|.
\end{equation}
\end{lemma}
\begin{proof}
Assume, without loss of generality, that $tr_{\otimes}(A)\geq tr_{\otimes}(B)$. Thus
\begin{eqnarray*}
tr_{\otimes}\vert A- B \vert=\displaystyle\max_{1\leq i \leq n}\vert a_{ii}-b_{ii} \vert&\geq& \vert a_{ii}-b_{ii} \vert\\
&\geq& a_{ii}-b_{ii}\ \ \ \ \mathrm{for~all}~i\in \{1, \ldots, n \}.
\end{eqnarray*}
Hence $tr_{\otimes}\vert A- B \vert+ b_{ii}\geq a_{ii}~ \mathrm{for~all}~i\in \{1, \ldots, n \}$ and so
\[
tr_{\otimes}\vert A- B \vert+ tr_{\otimes}(B)\geq tr_{\otimes}(A),
\]
which proves the first inequality in (\ref{Lip_of_tr}). The second inequality in (\ref{Lip_of_tr}) is trivial.
\end{proof}
The following lemma can be proved in a similar manner.
\begin{lemma}\label{lemma2}
If $A,B,C\in M_n(\mathbb{R}_+)$, then 
\[
\vert A\otimes B-A\otimes C\vert\leq A\otimes \vert B-C\vert .
\]
\[
\vert B\otimes A-C\otimes A\vert\leq\vert B-C\vert\otimes A.
\]
Consequently, 
\[
\Vert A\otimes B-A\otimes C\Vert\leq \Vert A\Vert \Vert B-C\Vert,
\]
and
\[
\Vert B\otimes A-C \otimes A\Vert\leq \Vert A\Vert \Vert B-C\Vert.
\]
\end{lemma}
The next fact is presumably known (see e.g. \cite[Lemma 1]{GPZ} for the case of square matrices). We include a proof here for the sake of completeness.  

\begin{lemma} \label{tr(AB=BA)} 
Let $A \in M_{n \times k}(\mathbb{R}_{+})$, $B \in M_{k \times n}(\mathbb{R}_{+})$. Then
$tr_{\otimes }(A\otimes B) = tr _{\otimes }(B \otimes A).$
 \end{lemma}
\begin{proof} 
 We have $(A\otimes B)_{ii} =  \max _{l = 1, \ldots , k} (a_{il}b_{li})$ for each $i=1, \ldots , n$ 
and \\
 $(B\otimes A)_{ll} =  \max _{i = 1, \ldots , n} (b_{li}a_{il})$ for each $l=1, \ldots , k$. 
Therefore $tr_{\otimes }(A\otimes B) =\max_i (A\otimes B)_{ii} =\max _l (B\otimes A)_{ll}= tr _{\otimes }(B \otimes A) $, which completes the proof.
\end{proof}


\begin{proposition}\label{InequalityLipschitzcontinuousA}
Let $A\in M_n(\mathbb{R}_+)$ and let $1\leq k \leq n$ be a positive integer. Consider the map  $f_{A}: \mathcal{X}_{n\times k}\longrightarrow \mathbb{R}_{+},$ where
\begin{eqnarray*}
 f_{A}(X):=tr_{\otimes}(X^t\otimes A\otimes X).
\end{eqnarray*}
Then
\begin{eqnarray}
\nonumber
\vert f_{A}(X)-f_{A}(Y) \vert & \leq &
  \Vert A \Vert \Vert X\otimes X^{t}-Y\otimes Y^{t} \Vert \\
 & \leq & 
  \|A\| (\|X\| + \|Y\|) \|X-Y\|
 \label{good_L}
\end{eqnarray}
for all $X, Y \in \mathcal{X}_{n\times k}$.
\end{proposition}
\begin{proof}
Let $X, Y\in \mathcal{X}_{n\times k}$. By  Lemmas \ref{lemma1}, \ref{tr(AB=BA)} and \ref{lemma2} we have
\begin{eqnarray*}
\vert f_{A}(X)-f_{A}(Y) \vert &=& \vert tr_{\otimes}(X^{t}\otimes A\otimes X)- tr_{\otimes}(Y^{t}\otimes A \otimes Y)\vert \\
&=& \vert tr_{\otimes}(X\otimes X^{t}\otimes A)- tr_{\otimes}(Y\otimes Y^{t}\otimes A)\vert\\
&\leq & \Vert X\otimes X^{t}\otimes A- Y\otimes Y^{t}\otimes A \Vert \\
&\le &
 \Vert  \vert X\otimes X^{t}- Y\otimes Y^{t}\vert\otimes A\Vert \\
&\leq & 
  \Vert A \Vert \Vert X\otimes X^{t}-Y\otimes Y^{t} \Vert,
\end{eqnarray*}
which proves the first inequality in (\ref{good_L}). Since
\begin{eqnarray*}
\Vert X\otimes X^{t}-Y\otimes Y^{t}\Vert &=& \Vert X\otimes X^{t}- Y\otimes X^{t}+Y\otimes X^{t}-Y\otimes Y^{t}\Vert \\
&\leq & \Vert X\otimes X^{t}-Y\otimes X^{t}\Vert+\Vert Y\otimes X^{t}-Y\otimes Y^{t}\Vert \\
&\leq & \Vert X\Vert\Vert X-Y\Vert+\Vert Y\Vert \Vert X-Y\Vert.
\end{eqnarray*}
Therefore
\begin{equation*}
\Vert X\otimes X^{t}-Y\otimes Y^{t}\Vert\leq\displaystyle \left(
\Vert X\Vert+\Vert Y\Vert
\right)\Vert X-Y\Vert,
\end{equation*} 
which verifies the second inequality in (\ref{good_L}).
\end{proof}
\begin{corollary} Let $A\in M_n(\mathbb{R}_+)$ and let $1\leq k \leq n$ be a positive integer. 
 For each $Z\in \mathcal{X}_{n\times k}$ and $X,Y\in \mathcal{X}_{n\times k}$ such that 
 $\|X-Z\|\le \frac{1}{2}$ and $\|Y-Z\|\le \frac{1}{2}$
  we have 
 \[
 \vert f_{A}(X)-f_{A}(Y) \vert\leq 
 \Vert A \Vert \displaystyle \left(
2\Vert Z\Vert+1
\right)\Vert X-Y\Vert.
 \] 
 Therefore $f_{A}: \mathcal{X}_{n\times k}\longrightarrow \mathbb{R}_{+}$,  is locally Lipschitz continuous. 
\end{corollary}
\begin{proof} 
From the assumed inequalities
$\Vert X-Z\Vert  \leq \frac{1}{2},~~ \Vert Y-Z\Vert  \leq \frac{1}{2}$
 it follows that 
$\Vert X\Vert \leq \Vert Z\Vert+\frac{1}{2}~~\textrm{and}~ \Vert Y\Vert \leq \Vert Z\Vert +\frac{1}{2}.$
This, together with (\ref{good_L}), proves the desired result.
\end{proof}
Our next aim is to rigorously prove  that equality (\ref{TO_CORRECT}) holds for all $1 \le k< n$ (Theorem \ref{Main Theorem}). We first establish the next two related results.
\begin{lemma}
 Suppose that $A=(a_{ij})\in M_n(\mathbb{R}_+)$ and let $1\leq k< n$ be a positive integer. We have 
 \[
 W_{max}^{k}(A)\subseteq [c, d] \;\;\mathit{and} \;\; \{c, d\} \subset W_{max}^{k}(A),
 \]
 where $c=\min\{\displaystyle\oplus_{j=1}^ka_{i_ji_j}:~1\leq i_1<i_2<\cdots <i_k\leq n\}$ and $d=\displaystyle\max_{1\leq i,j\leq n}a_{ij} =\|A\|$.  
 \label{main_correction2}
 \end{lemma}
 \begin{proof}
 Throughout the  proof we may by Theorem \ref{TheoremFirstknumer} (ii)
assume  that
  \begin{equation}
  \label{order}
  A=\left[
\begin{array}{cccc}
a_{11}&a_{12}&\cdots&a_{1n}\\
a_{21}&a_{22}&\cdots&a_{2n}\\
\vdots&\vdots&\ddots&\vdots\\
a_{n1}&a_{n2}&\cdots&a_{nn}
\end{array}
\right],
\end{equation}
where $a_{11}\leq a_{22}\leq \cdots \leq a_{nn}$.  For $k=1,$ we have $W_{max}^{1}(A)=W_{max}(A)=[c,  d],$ where $c=\displaystyle\min_{1\leq i\leq n}a_{ii},   d=\displaystyle\max_{1\leq i, j\leq n}a_{ij},$  \cite[Theorem 2]{TZAH}.  For $1 <k < n$ we first establish that $W_{max}^{k}(A)\subseteq [c, d]$.  Let $z \in W_{max}^{k}(A)$. 
Then
  \[
  z=\displaystyle\max_{1\leq r\leq k}(\max_{1\leq i, j \leq n}x_{i}^{(r)}x_{j}^{(r)}a_{ij})
  \]
  for some $X=[x^{(1)},\dots ,x^{(k)}]\in\mathcal{X}_{n\times k}$. By Lemma \ref{properties}   $\displaystyle \max\{x^{(r)}_{1}, x^{(r)}_{2}, \ldots, x^{(r)}_{n} \}=1$ for each   $1\leq r\leq k$ and   $x^{(r)}_{i}x^{(s)}_{i}=0$ for $ 1\leq r\neq s\leq k$ and all $1\leq i\leq n$.
Choose  $1\leq r_{1}\leq k, 1\leq i_{r_{1}}, j_{r_{1}}\leq n$ such that
  \[
  z=\displaystyle\max_{1\leq r\leq k}(\max_{1\leq i , j \leq n}x_{i}^{(r)}x_{j}^{(r)}a_{ij})=x_{i_{r_{1}}}^{(r_{1})}x_{j_{r_{1}}}^{(r_{1})}a_{i_{r_{1}}j_{r_{1}}}.
  \]
   Then by Lemma \ref{properties} (ii)
  \[
 z=x_{i_{r_{1}}}^{(r_{1})}x_{j_{r_{1}}}^{(r_{1})}a_{i_{r_{1}}j_{r_{1}}}\leq a_{i_{r_{1}}j_{r_{1}}}\leq d.
  \]
It is sufficient to show that $z\geq a_{kk}$. If $a_{kk} =0$ this inequality is obvious, so we may assume that $a_{kk}>0$.  For the sake of a contradiction, suppose that $z<a_{kk}$. Then
\[
x_{j}^{(i)}x_{j}^{(i)}a_{jj}\leq z<a_{jj}~~ \forall~ 1\leq i\leq k,~ \forall~ k\leq j\leq n.
\]
Since $a_{jj}>0$,  $ \forall~k\leq j\leq n,$ it follows that $x_{j}^{(i)}\neq 1$, $~ ~\forall~ 1\leq i\leq k,~ \forall~ k\leq j\leq n.$ This contradicts that fact that $X\in \mathcal{X}_{n\times k}$, by Lemma \ref{properties} (iii). Thus $z\geq a_{kk}$ and hence $W_{max}^{k}(A)\subseteq [c, d]$.

We conclude the proof by establishing that $\{c, d\} \subset W_{max}^{k}(A)$. 
Let
  $X=[x^{(1)}, x^{(2)}, \ldots, x^{(k)}], $
 where $x^{(i)}=e_{i},~~~ 1\leq i \leq k.$
 Then $X\in \mathcal{X}_{n\times k}$ and  $tr_{\otimes}(X^{t}\otimes A \otimes X)=a_{kk}=c \in W_{max}^{k}(A)$. Let $d=\displaystyle\max_{1\leq i,j\leq n}a_{ij}=a_{rs}$. 
    Let
 \[
  y^{(1)}=[0,  \ldots, 0, 1, 0, \ldots, 0, 1,0, \ldots, 0]^{t},~~ y_{r}^{(1)}=y_{s}^{(1)}=1,
 \]
and 
\begin{equation*}
Y=[y^{(1)},e_{i_1},\cdots ,e_{i_{k-1}}],
\end{equation*} 
such that $1\leq i_1<i_2<\cdots <i_{k-1}\leq n$ and $i_j\neq r,s$ for $j=1,\cdots ,k-1$.
Then 
  $Y\in \mathcal{X}_{n\times k}~\mathrm{and}~tr_{\otimes}(Y^{t}\otimes A \otimes Y)=a_{rs}=d \in W_{max}^{k}(A),$ 
  which completes the proof. 
 \end{proof}
\begin{lemma}\label{Total_Lemmas}
Suppose that  $A=(a_{ij})\in M_n(\mathbb{R}_+)$  and let $1\leq k< n$ be a positive integer. We have
 \[
 [c, d]\subseteq W_{max}^{k}(A),
 \]
 where $c=\min\{\displaystyle\oplus_{j=1}^ka_{i_ji_j}:~1\leq i_1<i_2<\cdots <i_k\leq n\}$ and $d=\displaystyle\max_{1\leq i,j\leq n}a_{ij} =\|A\|$.
\end{lemma}
\begin{proof} 
Throughout the  proof we may by Theorem \ref{TheoremFirstknumer} (ii)
assume  that
  \begin{equation}
  \label{order}
  A=\left[
\begin{array}{cccc}
a_{11}&a_{12}&\cdots&a_{1n}\\
a_{21}&a_{22}&\cdots&a_{2n}\\
\vdots&\vdots&\ddots&\vdots\\
a_{n1}&a_{n2}&\cdots&a_{nn}
\end{array}
\right],
\end{equation}
where $a_{11}\leq a_{22}\leq \cdots \leq a_{nn}$.  For $k=1,$ we have $W_{max}^{1}(A)=W_{max}(A)=[c,  d],$ where $c=\displaystyle\min_{1\leq i\leq n}a_{ii},   d=\displaystyle\max_{1\leq i, j\leq n}a_{ij},$  \cite[Theorem 2]{TZAH}.  So assume that $1 < k <n$ which leads to $c=a_{kk}$. Now,
let $z\in [a_{kk},a_{rs}]$ be fixed. By Theorem \ref{TheoremFirstknumer} (iv) assume that $r\le s$. We consider two cases.
\begin{itemize}
\item[Case 1:]  For $k<s,$ by letting $t=\max\{r,k\}$  we  distinguish two additional subcases.

(i) If $a_{tt}\leq z\leq a_{rs}$, then we define $X=[x^{(1)},\ldots ,x^{(k)}]$ with
\[
x^{(1)}=[0,\ldots ,1,\ldots ,0,\frac{z}{a_{rs}},0,\ldots ,0],\ x_r^{(1)}=1,x_s^{(1)}=\frac{z}{a_{rs}},
\]
\begin{equation*}
x^{(i)}=\begin{cases}
e_{i-1}&2\leq i\leq r\\
e_i&r<i\leq k.
\end{cases}
\end{equation*}
In this case it follows that 
\[
X\in \mathcal{X}_{n\times k},~ tr_{\otimes}(X^{t}\otimes A \otimes X)=z.
\]
(ii) If $a_{kk}\leq z< a_{tt},$  then we have two cases.

(a) If $\max\{a_{(k-1)t},a_{t(k-1)}\}\leq a_{kk}$, by taking $X=[x^{(1)},\ldots ,x^{(k)}]$, where 
\[
x^{(1)}=[0,\ldots ,1,\ldots ,0,\sqrt{\frac{z}{a_{tt}}},0,\ldots ,0],\ x_{k-1}^{(1)}=1,x_t^{(1)}=\sqrt{\frac{z}{a_{tt}}},
\]
\begin{equation*}
x^{(i)}=\begin{cases}
e_{i-1}&2\leq i\leq k-1\\
e_i&i=k
\end{cases}.
\end{equation*}
 We have
\[
X\in \mathcal{X}_{n\times k},~ tr_{\otimes}(X^{t}\otimes A \otimes X)=z.
\]
(b) If $p=\max\{a_{(k-1)t}, a_{t(k-1)}\}> a_{kk},$ then we have two further subcases.

$\bullet$ If $a_{kk}\leq z\leq \frac{p^{2}}{a_{tt}},$  then by taking $X=[x^{(1)},\ldots ,x^{(k)}]$, where 
\[
x^{(1)}=[0,\ldots ,1,\ldots ,0,\frac{z}{p},0,\ldots ,0],\ x_{k-1}^{(1)}=1,x_t^{(1)}=\frac{z}{p},
\]
\begin{equation*}
x^{(i)}=\begin{cases}
e_{i-1}&2\leq i\leq k-1\\
e_i&i=k
\end{cases}.
\end{equation*}

we have
\[
X\in \mathcal{X}_{n\times k},~ tr_{\otimes}(X^{t}\otimes A \otimes X)=z.
\]
$\bullet\bullet$ If $\frac{p^{2}}{a_{tt}}< z< a_{tt},$ then  by taking $X=[x^{(1)},\ldots ,x^{(k)}]$, where 
\[
x^{(1)}=[0,\ldots ,1,\ldots ,0,\sqrt{\frac{z}{a_{tt}}},0,\ldots ,0],\ x_{k-1}^{(1)}=1,x_t^{(1)}=\sqrt{\frac{z}{a_{tt}}},
\]
\begin{equation*}
x^{(i)}=\begin{cases}
e_{i-1}&2\leq i\leq k-1\\
e_i&i=k ,
\end{cases}.
\end{equation*}

We have
\[
X\in \mathcal{X}_{n\times k},~ tr_{\otimes}(X^{t}\otimes A \otimes X)=z.
\]
\item[Case 2:]If $k\geq s,$  then we have two related subcases.

$(i)^{\prime}$ If $a_{(k+1)(k+1)}\leq z\leq a_{rs},$ then by taking $X=[x^{(1)},\ldots ,x^{(k)}]$, where 
\[
x^{(1)}=[0,\ldots ,1,\ldots ,0,\ldots,\frac{z}{a_{rs}},0,\ldots ,0],\ x_{r}^{(1)}=1,x_s^{(1)}=\frac{z}{a_{rs}},
\]
\begin{equation*}
x^{(i)}=\begin{cases}
e_{i-1}&2\leq i\leq r\\
e_i&r<i<s\\
e_{i+1}&s\leq i\leq k,
\end{cases}
\end{equation*}
we have
\[
X\in \mathcal{X}_{n\times k},~ tr_{\otimes}(X^{t}\otimes A \otimes X)=z.
\]
$(ii)^{\prime}$ If $a_{kk}\leq z< a_{(k+1)(k+1)}$, 
 then we have two cases.

(a) If $\max\{a_{(k-1)(k+1)},a_{(k+1)(k-1)}\}\leq a_{kk}$, by taking $X=[x^{(1)},\ldots ,x^{(k)}]$, where 
\[
x^{(1)}=[0,\ldots ,1,0,\sqrt{\frac{z}{a_{(k+1)(k+1)}}},0,\ldots ,0],\ x_{k-1}^{(1)}=1, x_{k+1}^{(1)}=\sqrt{\frac{z}{a_{(k+1)(k+1)}}},
\]
\begin{equation*}
x^{(i)}=\begin{cases}
e_{i-1}&2\leq i\leq k-1\\
e_i&i=k,
\end{cases}
\end{equation*}
we have
\[
X\in \mathcal{X}_{n\times k},~ tr_{\otimes}(X^{t}\otimes A \otimes X)=z.
\]
(b) If $p=\max\{a_{(k-1)(k+1)}, a_{(k+1)(k-1)}\}> a_{kk},$ then we consider further two subcases.

$\bullet$ If $a_{kk}\leq z\leq \frac{p^{2}}{a_{(k+1)(k+1)}},$ then  by taking $X=[x^{(1)},\ldots ,x^{(k)}]$, where 
\[
x^{(1)}=[0,\ldots ,1 ,0,\frac{z}{p},0,\ldots ,0],\ x_{k-1}^{(1)}=1,x_{k+1}^{(1)}=\frac{z}{p},
\]
\begin{equation*}
x^{(i)}=\begin{cases}
e_{i-1}&2\leq i\leq k-1\\
e_i&i=k
\end{cases}
\end{equation*}
we have
\[
X\in \mathcal{X}_{n\times k},~ tr_{\otimes}(X^{t}\otimes A \otimes X)=z.
\]
$\bullet\bullet$ If $\frac{p^{2}}{a_{(k+1)(k+1)}}< z< a_{(k+1)(k+1)},$  by taking $X=[x^{(1)},\ldots ,x^{(k)}]$, where 
\[
x^{(1)}=[0,\ldots ,1,\ldots ,0,\sqrt{\frac{z}{a_{(k+1)(k+1)}}},0,\ldots ,0],\ x_{k-1}^{(1)}=1,x_{k+1}^{(1)}=\sqrt{\frac{z}{a_{(k+1)(k+1)}}},
\]
\begin{equation*}
x^{(i)}=\begin{cases}
e_{i-1}&2\leq i\leq k-1\\
e_i&i=k,
\end{cases}
\end{equation*}
we have
\[
X\in \mathcal{X}_{n\times k},~ tr_{\otimes}(X^{t}\otimes A \otimes X)=z,
\]
\end{itemize}
which completes the proof.
\end{proof}
 \begin{theorem}\label{Main Theorem}
 Suppose that $A=(a_{ij})\in M_n(\mathbb{R}_+)$ and let $1\leq k< n$ be a positive integer. Then
 \[
 W_{\max}^k(A)=[c,d],~~1\leq k<n,
 \]
 where $c=\min\{\displaystyle\oplus_{j=1}^ka_{i_ji_j}:~1\leq i_1<i_2<\cdots <i_k\leq n\}$ and $d=\displaystyle\max_{1\leq i,j\leq n}a_{ij}$.  Moreover, 
 \[
 W_{\max}^n(A)=\{\displaystyle\max_{1\leq i\leq n}a_{ii} \}.
 \]
 \label{main_correction}
 \end{theorem}
 \begin{proof}
 The result follows from Theorem \ref{wmaxA},   Lemmas \ref{main_correction2} and \ref{Total_Lemmas} and from (\ref{WmaxNnumberone}).
 \end{proof}
\begin{remark}\label{Exaple3TZAH}
In \cite[Example 3]{TZAH},  for the matrix 
\[
A=\left [
\begin{array}{cccc}
4 & 7 & 5 & 8\\
8 & 2 & 0 & 7\\
2 & 8 & 1 & 4\\
1 & 6 & 2 & 2
\end{array}
\right].
\]
%
it was stated that $W_{max}^{4}(A)=[4,  8],$ which is  of course not correct. From (\ref{WmaxNnumberone}) it follows that $W_{max}^{4}(A)=\{4\}$.
\end{remark}
\begin{example}
Let $A=(a_{ij})\in M_{n\times n}(\mathbb{R}_{+}),$ where $a_{11}\leq a_{22}\leq \cdots \leq a_{nn}$ and let $\displaystyle\max_{1\leq i, j\leq n}a_{ij}=a_{rs}$. By Theorem \ref{Main Theorem}, we have $W_{max}^{k}(A)=[a_{kk},  a_{rs}],$  where $1\leq k<n$ and $W_{max}^{n}(A)=\{a_{nn}\}$.

\end{example}
\begin{example}
Let
\[
A=\left [
\begin{array}{cccccc}
2.5&5.2&4.1&2.3&4&3.5\\
5&3&6.2&3&3.5&4.7\\
3.7&4&5.2&6&5.8& 4.3\\
2.5&6&1.7&6.2&9& 8.1\\
7.2&5.3&4.2&6.1&7.4 &7\\
8.1& 7.6& 5.9& 3.8& 9& 8.3
\end{array}
\right].
\]
Then $\displaystyle \max_{1\leq i, j\leq 6}a_{ij}=a_{45}=9$. By Theorem \ref{Main Theorem} we have
\[
W_{max}^{1}(A)=W_{max}(A)=[2.5,  9],~ W_{max}^{2}(A)=[3,  9],~W_{max}^{3}(A)=[5.2,  9],
\]
\[ 
W_{max}^{4}(A)=[6.2,  9],~W_{max}^{5}(A)=[7.4,  9] ~ \textit{and}~ W_{max}^{6}(A)=\{8.3\}.
\]

\end{example}
\bigskip

Let $A\in M_{n\times n}(\mathbb{R}_{+})$ and $x\in\mathbb{R}_+^n$. Then we 
let $r_{x}(A)$ denote the local spectral radius of $A$ at $x$, i.e., $r_{x}(A)=\displaystyle\limsup_{j\rightarrow \infty} \Vert A_{\otimes}^{j}\otimes x \Vert^{\frac{1}{j}}$. It was shown in \cite{M.PEPERKO1} for $ x=(x_{1}, \ldots, x_{n})^{t}\in \mathbb{R}_{+}^{n}, x\neq 0$ it holds $r_{x}(A)=\displaystyle\lim_{j\rightarrow\infty}\Vert A_{\otimes}^{j}\otimes x \Vert^{\frac{1}{j}}$ and  $r_{x}(A)=\max\{r_{e_{i}}(A): i=1, \ldots, n, x_{i}\neq 0\},$  
where $x_{i}$ denotes the $i$th coordinate of $x$. We say that  $\mu\geq 0$ is a geometric max eigenvalue of $A$ if $A\otimes  x=\mu x$ for some $x\neq 0$ with $x\geq 0$. Let $\sigma_{max}(A)$ denote the set of geometric max eigenvalues of $A$. It is known (see e.g \cite{Gu94}, \cite[Theorem 2.7]{M.PEPERKO1}) that 
\[
\sigma_{max}(A)=\{\mu: \mu=r_{e_{j}}(A), ~\textrm{for} ~ j\in \{1, \ldots, n\}\}.
\]
We define the standard vector multiplicity of geometric max eigenvalue $\mu$ as the number of indices $j$ such that $\mu=r_{e_{j}}(A)$.

The role of the spectral radius of $A$ in max algebra is played by the maximum cycle geometric mean $\mu(A),$ which is defined by
\begin{equation}\label{definition_mu(A)}
\mu(A)=\displaystyle\max\left\lbrace (a_{i_{1}i_{k}}\ldots a_{i_{3}i_{2}}a_{i_{2}i_{1}})^{\frac{1}{k}}: k\in \mathbb{N }~\textit{and}~ i_{1}, \ldots, i_{k}\in \{1, \ldots, n\} \right\rbrace,
\end{equation}
and is equal to 
\begin{equation*}
\mu(A)=\max\left\lbrace (a_{i_{1}i_{k}}\ldots a_{i_{3}i_{2}}a_{i_{2}i_{1}})^{\frac{1}{k}}: k\leq n~ \textit{and distinct}~i_{1}, \ldots, i_{k}\in \{1, \ldots, n\}\right\rbrace.
\end{equation*}
It is known (see \cite{Butkovic})   that $\mu(A)$ is the largest geometric max eigenvalue of $A,$  i.e.,
$\mu(A)=\max\{\mu: \mu\in \sigma_{max}(A)\}$ and thus $\mu(A)=\displaystyle\max_{j=1, \ldots, n}r_{e_{j}}(A)$.  

The max permanent of $A$  is given by
\[
perm(A)=\displaystyle\max_{\sigma\in \sigma_{n}}a_{1\sigma(1)}\ldots a_{n\sigma(n)},
\] 
where $\sigma_{n}$ is the group of permutations on $\{1, \ldots, n\}$. The characteristic maxpolynomial of $A$ (see e.g. \cite{Butkovic,  RosenmannLehnerPeperko, TS})  is a max polynomial
\[
\mathcal{X}_{A}(x)=perm(xI\oplus A).
\]
Clearly $\mathcal{X}_{P^T\otimes A\otimes P}(x)=\mathcal{X}_{A}(x)$ holds for each $n\times n$ permutation matrix $P$ and $x\ge 0$.
We call its tropical roots (the points of nondifferentiability of $\mathcal{X}_{A}(x)$ considered as a function on $[0, \infty)$) the algebraic max eigenvalues (or also tropical eigenvalues) of $A$. The set of all algebraic max eigenvalues is denoted by $\sigma_{trop}(A)$.  For $\lambda\in \sigma_{trop}(A)$ its multiplicity, as a tropical root of $\mathcal{X}_{A}(x)$ (see e.g \cite{Butkovic,  RosenmannLehnerPeperko, TS}),
is called the algebraic multiplicity of $\lambda$.  It is well known that $\sigma_{max}(A)\subset \sigma_{trop}(A)$\cite[Remark 2.3]{TS}  and that $\mu(A)=\max\{\lambda: \lambda\in \sigma_{trop}(A)\}.$  However, in general, the sets $\sigma_{max}(A)$ and $\sigma_{trop}(A)$ may not coincide. Further, 
\begin{equation}
\label{p_eq}
\sigma_{max}(P^t \otimes A\otimes P)=\sigma_{max}(A)\;\;\mathrm{and}\;\; \sigma_{trop}(P^t \otimes A\otimes P)=\sigma_{trop}(A)
\end{equation}
for any  $n\times n$ permutation matrix $P$.

Recall that the max convex hull of a set $M\subseteq \mathbb{R}_{+},$ which is denoted by $conv_{\otimes}(M),$ is defined as the set of all max convex linear combinations of elements from $M,$ i.e., 
\[
conv_{\otimes}(M):= \left\lbrace \displaystyle\bigoplus_{i=1}^{m}\alpha_{i}x_{i}: m\in \mathbb{N},  x_{i}\in M, \alpha_{i}\geq 0, i=1, \ldots, m, \bigoplus_{i=1}^{m}\alpha_{i}=1\right\rbrace.
\] 

In \cite[Section 4]{TZAH},   we defined the max $k-$geometric spectrum and $k-$tropical spectrum of $A \in M_{n}(\mathbb{R}_{+})$ as follows.
Let $A \in M_{n}(\mathbb{R}_{+}), 1\leq k \le  n$, let $\mu_{1}, ..., \mu_{n} \in \sigma_{max}(A)$ counting standard vector multiplicities and  let $\lambda_{1}, ..., \lambda_{n} \in \sigma_{trop}(A)$ counting tropical multiplicities.  The max $k-$geometric spectrum of $A$ is given by
\begin{equation*}
\sigma_{max}^{k}(A)=\displaystyle\left\lbrace\displaystyle\bigoplus_{j=1}^{k}\mu_{i_{j}}: 1\leq i_{1}<i_{2}<\cdots <i_{k}\leq n \right\rbrace.
\end{equation*}
and the  $k-$tropical max spectrum of $A$ is 
\begin{equation*}
\sigma_{trop}^{k}(A)=\displaystyle\left\lbrace\displaystyle\bigoplus_{j=1}^{k}\lambda_{i_{j}}: 1\leq i_{1}<i_{2}<\cdots <i_{k}\leq n \right\rbrace.
\end{equation*}
It is clear that $\sigma_{max}^{1}(A)=\sigma_{max}(A)$  and  $\sigma_{trop}^{1}(A)=\sigma_{trop}(A)$. By (\ref{p_eq}) and from the definitions above we have 
\begin{equation}
\label{p_eq2}
\sigma_{max} ^{k} (P^t \otimes A\otimes P)=\sigma_{max} ^{k}(A)\;\;\mathrm{and}\;\; \sigma_{trop} ^{k}(P^t \otimes A\otimes P)=\sigma_{trop} ^{k}(A),
\end{equation}
for any  $n\times n$ permutation matrix $P$ and all $k=1, \ldots , n$.
\begin{remark}\label{proposition7TZH}
{\rm 
It was stated in \cite[Proposition 7]{TZAH} that 
$conv_{\otimes}(\sigma_{max}^{n}(A))\subseteq W_{max}^{n}(A)$  and $conv_{\otimes}(\sigma_{trop}^{n}(A))\subseteq W_{max}^{n}(A)$. 
However, this turns out to be false in general. Indeed, from the definitions above we have
\[
conv_{\otimes}(\sigma_{max}^{n}(A))=conv_{\otimes}(\sigma_{trop}^{n}(A))=\{\mu(A)\}.
\]
On the other hand, for \[
A=\left [
\begin{array}{cc}
0&1\\
1&0
\end{array}
\right]\in M_{2\times 2}(\mathbb{R}_{+}).
\]
we have $\mu(A)=1,$  and
\[
conv_{\otimes}(\sigma_{max}^{2}(A))=conv_{\otimes}(\sigma_{trop}^{2}(A))=\{1\}\nsubseteq W_{max}^{2}(A)=\{0\}.
\]
}
\label{conv_cor}
\end{remark}
However, we rigorously prove in Proposition \ref{inclusionsigmak} and Theorem \ref{trop_cor} below that the inclusions  $conv_{\otimes}(\sigma_{max}^{k}(A))\subseteq W_{max}^{k}(A)$  and $conv_{\otimes}(\sigma_{trop}^{k}(A))\subseteq W_{max}^{k}(A)$ hold for all  $1 \le k <n$ (as stated in  \cite[Proposition 7]{TZAH}). For this we need to recall some related facts.

Every nonnegative matrix $A=(a_{ij})\in M_{n}(\mathbb{R}_{+})$ can be transformed in linear time by simultaneous permutations of the rows and columns into its \textit{Frobenius normal form} (\textbf{FNF}) \cite{BapatRaghavan, BrualdiRyser, Rosen}
\begin{equation}\label{FNF}
\left [
\begin{array}{ccccc}
B_{l}&0 & 0 & \ldots & 0\\
\star &B_{l-1} & 0 &  \ldots & 0\\
 \vdots & \vdots & \ddots & \vdots & 0\\
\star&\star & \star & \ldots & B_{1}
\end{array}
\right],
\end{equation}
where $B_{1}, \ldots, B_{l}$ are irreducible square submatrices of $A$ or $1 \times 1$ zero blocks. The diagonal blocks are determined uniquely up to a simultaneous permutation of their rows and columns: however, their order is not determined uniquely.  Observe that the \textbf{FNF} is a particularly convenient form for studying certain spectral properties of nonnegative matrices. Since these are essentially preserved by simultaneous permutations of the rows and columns (\cite[Proposition 4.1.3]{Butkovic}, (\ref{p_eq2}))
we will often assume, without loss of generality, that the matrix under consideration is in \textbf{FNF}.

If $A$ is in \textbf{FNF}, then the corresponding partition of the node set $N$ of $\mathcal{R}(A)$ is denoted as $N_{1}, \ldots, N_{l}$, partitioned conformally with the subsets corresponding to the irreducible submatrices in the \textbf{FNF} of $A$. 
Then the induced subgraphs $\mathcal{R}(A)[N_{\mu}] (\mu=1, \ldots, l)$ are strongly connected and an arc from $N_{\mu}$ to $N_{\nu}$ in $\mathcal{R}(A)$ exists only if $\mu\leq \nu$ (see also \cite{BrualdiRyser}).

Clearly, every $B_{\mu}$ has a unique max geometric eigenvalue $r_{\otimes}(B_{\mu})$. As a slight abuse of language we will, for simplicitly, also say that $r_{\otimes}(B_{\mu})$ is the max geometric  eigenvalue of $N_{\mu}$. 

The reduced graph denoted by $\mathcal{R}(A)$ is a digraph whose nodes correspond to $N_{\mu}$ for $\mu=1, \ldots, l$ and the set of edges is 
\[
\{(\mu, \nu): \mathrm{there~ exist}~k\in N_{\mu}~\mathrm{and}~ j\in N_{\nu} ~\mathrm{such~that}~ a_{kj}>0 \}.
\]  
By a class of $A$ we mean a node $\mu$  (or also the corresponding set $N_{\mu}$) of the reduced graph $\mathcal{R}(A)$. A class $\mu$ is trivial if $B_{\mu}$ is the $1\times 1$ zero matrix. Class $\mu$ accesses class $\nu,$  denoted $\mu\rightarrow \nu,$  if $\mu=\nu$  or if there exists a $\mu - \nu$ path in $\mathcal{R}(A)$ (a path that starts in $\mu$ and ends in $\nu$). A node $j$ of $\mathcal{R}(A)$ is accessed by a class $\mu,$ denoted by $\mu\rightarrow j,$ if $j$ belongs to a class $\nu$ such that $\mu\rightarrow \nu$.
The max eigenvalues $r_{e_{j}}(A)$ are described in the following way via this access relation, which is also an equivalence relation (see e.g. \cite[Corollary 2.9]{M.PEPERKO1}):
\begin{equation}
r_{e_{j}}(A)=\max\{r_{\otimes}(B_{\mu}): \mu\rightarrow j\}
\label{access}
\end{equation}
for all $j=1, \ldots, n$. 
For each $j=1, \ldots, n$ we have $r_{e_{j}}(A)=r_{\otimes}(B_{\nu})$ for some class $\nu$. The converse of this statement need not hold in general. Similarly, (see e.g. \cite[Corollary 2.10]{M.PEPERKO1})
\[
\sigma_{max}(A)=\displaystyle\left\lbrace r_{\otimes}(B_{\nu}): r_{\otimes}(B_{\nu})=\max \{r_{\otimes}(B_{\mu}): \mu\rightarrow \nu\}\right\rbrace.
\]
\begin{proposition}\label{inclusionsigmak}
Let $A\in M_{n}(\mathbb{R}_{+})$  and $1\leq k<n$. Then $\displaystyle conv_{\otimes}\left(\sigma_{max}^{k}(A)\right) \subseteq W_{max}^{k}(A)$.  
\end{proposition}
\begin{proof}
Without loss of generality, assume that
 $A=(a_{ij})$ is in a  \textbf{FNF} (\ref{FNF}). Assume that 
$r_{e_{i_1}}(A)\leq r_{e_{i_2}}(A)\leq \ldots \leq r_{e_{i_n}}(A)$ are the max geometric eigenvalues of $A$ counting standard vector multiplicities.  So 
\[
\sigma_{max}^{k}(A)=\{r_{e_{i_k}}(A), r_{e_{i_{k+1}}}(A), \ldots, r_{e_{i_n}}(A)\}.
\]
 For all $1\leq j\leq n,$  we have $r_{e_{i_j}}(A)=r_{\otimes}(B_{\nu})$ by (\ref{access}) for some class $\nu,  1\leq \nu \leq l$ such that $ \nu \rightarrow j$. Let $1\leq t\leq l$ be the smallest number such that 
 \[
 \{1, 2, \ldots, j\}\subseteq \displaystyle\bigcup_{\nu\in \{1, 2, \ldots, t\}}N_{\nu}.
 \]
   Hence for all $ j\in \{k, k+1, \ldots, n\}$ we have
\[
c\leq \displaystyle\bigoplus_{\nu=1}^{t}\bigoplus_{i=1}^{|N_{\nu}|}(B_{\nu})_{ii}\leq\displaystyle\bigoplus_{\nu=1}^{t}r_{\otimes}(B_{\nu})\leq \displaystyle\bigoplus_{s=1}^{j}r_{e_{i_s}}(A)=r_{e_{i_j}}(A)\leq d,
\]
 where 
 $c=\min\{\displaystyle\oplus_{j=1}^ka_{i_ji_j}:~1\leq i_1<i_2<\cdots <i_k\leq n\}$ and $d=\displaystyle\max_{1\leq i,j\leq n}a_{ij}.$  Thus
 $\sigma_{max}^{k}(A)\subseteq W_{max}^{k}(A)$. On the other hand  $W_{max}^{k}(A)$ is a closed interval and hence $\displaystyle conv_{\otimes}\left(  \sigma_{max}^{k}(A)\right) \subseteq W_{max}^{k}(A)$.
 
\end{proof}

For $A\in M_{n}(\mathbb{R}_{+})$ recall that the characteristic maxpolynomial $\mathcal{X}_{A}(x)$ equals $perm(xI\oplus A)$, where $x$ is an indeterminate, and can be written as
\begin{eqnarray*}
\mathcal{X}_{A}(x)&=&x^n\oplus\delta_1x^{n-1}\oplus\delta_2x^{n-2}\oplus\cdots\delta_{n-1}x\oplus\delta_n\\
&=&(x\oplus\lambda_1) (x\oplus\lambda_2)\cdots (x\oplus\lambda_n),
\end{eqnarray*}
where $\lambda_{1}\geq \lambda_{2}\geq \cdots \geq \lambda_{n}$ (\cite[Lemma 5.0.1]{Butkovic}) are the max algebraic eigenvalues (counting possible multiplicities)
{\color{black}
and where the coefficients $\delta_k$ 
can be chosen as  $\delta_{0}=1$ and
\[
\delta_k=\displaystyle\max_{B\in P_k(A)}perm(B),
\]
for $k=1,\cdots ,n$ (\cite[Theorem 5.3.2]{Butkovic}), where $P_k(A)$ is the set of all principal submatrices of $A$ of order $k$. 
}
For a given characteristic maxpolynomial, $\mathcal{X}_{A}(x)$ as above, we say that a coefficient $\delta_i$ is {\em inessential} (see also \cite{Butkovic}) if $\delta_i x^{n-i} \leq \sum_{j\neq i} \delta_j x^{n-j}$ for all $x \geq 0$; otherwise it is called {\em essential}.

 Let $A\in M_{n}(\mathbb{R}_{+}),$ where $a_{11}\leq a_{22}\leq \cdots \leq a_{nn}$. Then the characteristic maxpolynomial for $A$ can be written as
 \begin{eqnarray*}
\mathcal{X}_{A}(x)&=&x^n\oplus\delta_1x^{n-1}\oplus\delta_2x^{n-2}\oplus\cdots\delta_{n-1}x\oplus\delta_n\\
&=&x^{n}\oplus \delta_{i_{1}}x^{n-i_{1}}\oplus \delta_{i_{2}}x^{n-i_{2}} \oplus \cdots \oplus\delta_{i_{t}}x^{n-i_{t}},
\end{eqnarray*}
where $1 \leq t \leq n$, $1\leq i_1<i_2<\cdots <i_t\leq n$, and each of coefficients $\delta_{i_j}$, $j=1, \ldots , t$ are essential. In this case, it is known from  \cite{Butkovic} that the distinct nonzero max algebraic eigenvalues are given by
 \begin{eqnarray*}
 \lambda_{i_{r}}=
 \displaystyle(\frac{\delta_{i_{r}}}{\delta_{i_{r-1}}})^{\frac{1}{i_{r}-i_{r-1}}},~~~~~\mathit{ r=1, 2, \ldots, t, ~ i_{0}=0.}
 \end{eqnarray*}
 Furthermore, if $i_t < n$, then there is an additional max algebraic eigenvalue equal to 0 with corresponding multiplicity $n-i_t$. Thus the
 characteristic maxpolynomial for $A$ can also be written as
 \[ \mathcal{X}_{A}(x)= (x\oplus \lambda_{i_{1}})^{i_{1}}(x\oplus \lambda_{i_{2}})^{i_{2}-i_{1}}\cdots (x \oplus \lambda_{n})^{i_t-i_{t-1}}x^{n-i_t}.\]

\begin{remark}\label{r19}
Let $A\in M_n(\mathbb{R}_+)$ and let $1\leq i\leq n$ be a positive integer. Then the first and the last nonzero terms of $\mathcal{X}_{A}(x)$ are essential. 
\end{remark}
 
The following lemma is known (\cite[Lemma 5.0.1]{Butkovic}).
\begin{lemma}\label{orderroots}
Let $A\in M_n(\mathbb{R}_+)$ and let
\[
 \mathcal{X}_{A}(x)= x^{n}\oplus \delta_{i_{1}}x^{n-i_{1}}\oplus \delta_{i_{2}}x^{n-i_{2}} \oplus \cdots \oplus\delta_{i_{t}}x^{n-i_{t}},
\]
where $1\leq i_1<i_2<\cdots <i_t\leq n$ and let $\delta_{i_s} $ be  essential in $\mathcal{X}_{A}(x)$. Finally, let $i_0=0,~\delta_{i_0}=1.$
Then we have
\begin{equation}
\displaystyle\left(\frac{\delta_{i_{s}}}{\delta_{i_{s-1}}}\right)^{\frac{1}{i_{s}-i_{s-1}}}>\displaystyle\left(\frac{\delta_{i_{s+1}}}{\delta_{i_{s}}}\right)^{\frac{1}{i_{s+1}-i_{s}}}~~\mathit{for~ all}~ 1\leq s\leq t-1.
\end{equation}
\end{lemma}
{
\color{black}

\begin{remark}\label{peperko.concern}
Let $A\in M_n(\mathbb{R}_+)$ and let $1\leq i\leq n$ be a positive integer. If  $\delta_i=0$,  then there are at least $n-(i-1)$
 zeros on the main diagonal of $A$. Assume that $a_{11}\leq\cdots\leq a_{nn}$. So if $\delta_{i}=0$ then  $a_{n-i+1,n-i+1}=0$ and so $\delta_{i}=\delta_{i-1}a_{n-i+1,n-i+1}=0$.) 
\end{remark}
}
By (\ref{p_eq}), we may restrict to the case  $a_{11}\leq\cdots\leq a_{nn}$ in the following Lemma \ref{l1}, Remark \ref{specialorderingeig} and Theorem \ref{th4} below.

\begin{lemma}\label{l1}
Let $A\in M_{n}(\mathbb{R}_{+})$ such that $a_{11}\leq\cdots\leq a_{nn}$  and let $0\leq i< n$ be an  integer. Then  we have
\[
\delta_{i+1}\geq\delta_i a_{n-i,n-i}.
\]
\end{lemma}
\begin{proof}
In the case $i=0$ it is clear that $\delta_1=a_{nn}=\delta_0a_{nn}$.
Now let $1\leq i< n$ be given 
and  let $B\in M_{i}(\mathbb{R}_{+})$ be a principal submatrix of $A$ containing $\{j_{1},\ldots , j_{i}\}$ rows and columns of $A$ and $\delta_{i}=perm(B)$. By putting $t=\max (\{1, 2, \cdots, n\}\setminus \{j_{1}, j_{2}, \cdots, j_{i}\}) $, we have $t\geq n-i$ and
\[
\delta_{i+1}\geq a_{tt}\delta_i\geq a_{n-i,n-i}\delta_{i}.
\]
\end{proof}


\begin{remark}\label{specialorderingeig}
Suppose $A\in M_{n}(\mathbb{R}_{+}),$ where $a_{11}\leq a_{22}\leq \cdots \leq a_{nn}$. If all terms in the characteristic maxpolynomial are essential, then we have 
\begin{eqnarray*}
\mathcal{X}_{A}(x)&=&x^n\oplus\delta_1x^{n-1}\oplus\delta_2x^{n-2}\oplus\cdots\delta_{n-1}x\oplus\delta_n\\
&=&(x\oplus\lambda_1) (x\oplus\lambda_2)\cdots (x\oplus\lambda_n).
\end{eqnarray*}
In this case
\[
\lambda_i=\frac{\delta_i}{\delta_{i-1}},~i=1,\ldots ,n.
\]

Using Lemma \ref{l1} we have $\lambda_i\geq a_{n-i+1,n-i+1}$. Therefore in this special case, we have $\displaystyle conv_{\otimes}\left(\sigma_{trop}^{k}(A)\right)\subseteq W_{max}^{k}(A)$ for all $1\leq k< n$ by Theorem \ref{Main Theorem}.
\end{remark}


\begin{lemma}\label{lem1}
Let $A\in M_n(\mathbb{R}_+)$ and let
\[
 \mathcal{X}_{A}(x)= x^{n}\oplus \delta_{i_{1}}x^{n-i_{1}}\oplus \delta_{i_{2}}x^{n-i_{2}} \oplus \cdots \oplus\delta_{i_{t}}x^{n-i_{t}},
\]
where $1\leq i_1<i_2<\cdots <i_t\leq n$ and where $\delta_{i_j},~1\leq j\leq t$ are essential terms. Finally, let $i_0=0,~\delta_{i_0}=1,~\delta_{i_{t+1}}=0$.
Then we have
\begin{equation}\label{equ1}
\min\Bigg\{(\delta_{i_m})^{\frac{1}{i_m}},(\frac{\delta_{i_m}}{\delta_{i_1}})^{\frac{1}{i_{m}-i_1}},(\frac{\delta_{i_m}}{\delta_{i_2}})^{\frac{1}{i_{m}-i_2}},\cdots ,(\frac{\delta_{i_m}}{\delta_{i_{m-1}}})^{\frac{1}{i_{m}-i_{m-1}}}\Bigg \}=(\frac{\delta_{i_m}}{\delta_{i_{m-1}}})^{\frac{1}{i_{m}-i_{m-1}}},~0<m\leq t
\end{equation}
\begin{equation}\label{equ2}
\max\Bigg\{(\frac{\delta_{i_{m+1}}}{\delta_{i_{m}}})^{\frac{1}{i_{m+1}-i_{m}}},(\frac{\delta_{i_{m+2}}}{\delta_{i_{m}}})^{\frac{1}{i_{m+2}-i_{m}}},\cdots ,(\frac{\delta_{i_t}}{\delta_{i_{m}}})^{\frac{1}{i_{t}-i_{m}}}\Bigg \}=(\frac{\delta_{i_{m+1}}}{\delta_{i_{m}}})^{\frac{1}{i_{m+1}-i_{m}}},~0\leq m<t
\end{equation}
where $m$ is an integer.
\end{lemma}

\begin{proof}
For the first equation, let $0<m\leq t$ be given. If we assume that (\ref{equ1}) does not hold, then there is  $0\leq l<m-1$ such that 
\begin{equation}\label{equ3}
(\frac{\delta_{i_{m}}}{\delta_{i_{l}}})^{\frac{1}{i_{m}-i_{l}}}<(\frac{\delta_{i_{m}}}{\delta_{i_{m-1}}})^{\frac{1}{i_{m}-i_{m-1}}}.
\end{equation}
Since $\delta_{i_{m-1}}$ is essential, there is  $x$ such that
\[
\delta_{i_{l}}x^{n-i_{l}}\leq\delta_{i_{m-1}}x^{n-i_{m-1}},~\delta_{i_{m}}x^{n-i_{m}}\leq\delta_{i_{m-1}}x^{n-i_{m-1}}.
\]
So we have
\[
(\frac{\delta_{i_{m}}}{\delta_{i_{m-1}}})^{\frac{1}{i_{m}-i_{m-1}}}\leq x \le (\frac{\delta_{i_{m-1}}}{\delta_{i_{l}}})^{\frac{1}{i_{m-1}-i_{l}}},
\]
which contradicts (\ref{equ3}).

For the proof of (\ref{equ2}) let
$0\leq m<t$ be given. If we assume that (\ref{equ2}) does not hold, then there is  $m+1<l\leq t$ such that
\begin{equation}\label{equ4}
(\frac{\delta_{i_{m+1}}}{\delta_{i_{m}}})^{\frac{1}{i_{m+1}-i_{m}}}<(\frac{\delta_{i_{l}}}{\delta_{i_{m}}})^{\frac{1}{i_{l}-i_{m}}}.
\end{equation}
Since $\delta_{i_{m+1}}$ is essential, there is  $x$ such that
\[
\delta_{i_{l}}x^{n-i_{l}}\leq\delta_{i_{m+1}}x^{n-i_{m+1}},~\delta_{i_{m}}x^{n-i_{m}}\leq\delta_{i_{m+1}}x^{n-i_{m+1}}.
\]
So we have
\[
(\frac{\delta_{i_{l}}}{\delta_{i_{m+1}}})^{\frac{1}{i_{l}-i_{m+1}}}\leq x \le (\frac{\delta_{i_{m+1}}}{\delta_{i_{m}}})^{\frac{1}{i_{m+1}-i_{m}}},
\]
which contradicts (\ref{equ4}).
\end{proof}

\begin{lemma}\label{lem2}
Let $A\in M_n(\mathbb{R}_+)$ and let
\[
 \mathcal{X}_{A}(x)= x^{n}\oplus \delta_{i_{1}}x^{n-i_{1}}\oplus \delta_{i_{2}}x^{n-i_{2}} \oplus \cdots \oplus\delta_{i_{t}}x^{n-i_{t}},
\]
where $1\leq i_1<i_2<\cdots <i_t\leq n$ and where $\delta_{i_j},~1\leq j\leq t$ are  the essential terms. 
Moreover, let $0\leq m\leq t$ be an  integer, $\delta_{i_{-1}}=0$, $i_{-1}=-1$,   $i_0=0$  and also let $\delta_{i_0}=1$. Then for all $(\frac{\delta_{i_{m+1}}}{\delta_{i_{m}}})^{\frac{1}{i_{m+1}-i_{m}}}\leq x\leq (\frac{\delta_{i_{m}}}{\delta_{i_{m-1}}})^{\frac{1}{i_{m}-i_{m-1}}}$, we have
\[
x^n\oplus \delta_{i_1}x^{n-i_1}\oplus\delta_{i_2}x^{n-i_2}\oplus\cdots\oplus\delta_{i_t}x^{n-i_t}=\delta_{i_m}x^{n-i_m}.
\]
\end{lemma}

\begin{proof}
If the conclusion fails to hold, then there is $0\leq l\leq t,~l\neq m$ with $(\frac{\delta_{i_{m+1}}}{\delta_{i_{m}}})^{\frac{1}{i_{m+1}-i_{m}}}\leq x\leq (\frac{\delta_{i_{m}}}{\delta_{i_{m-1}}})^{\frac{1}{i_{m}-i_{m-1}}}$ such that
\[
\delta_{i_m}x^{n-i_m}<\delta_{i_l}x^{n-i_l}.
\]
Now, we have two cases.

Case (1): If $l>m$, then $x<(\frac{\delta_{i_{l}}}{\delta_{i_{m}}})^{\frac{1}{i_{l}-i_{m}}}$. By Lemma \ref{lem1} we have
\[
\max\Bigg\{(\frac{\delta_{i_{m+1}}}{\delta_{i_{m}}})^{\frac{1}{i_{m+1}-i_{m}}},(\frac{\delta_{i_{m+2}}}{\delta_{i_{m}}})^{\frac{1}{i_{m+2}-i_{m}}},\cdots ,(\frac{\delta_{i_t}}{\delta_{i_{m}}})^{\frac{1}{i_{t}-i_{m}}}\Bigg \}=(\frac{\delta_{i_{m+1}}}{\delta_{i_{m}}})^{\frac{1}{i_{m+1}-i_{m}}}.
\]
So $x<(\frac{\delta_{i_{m+1}}}{\delta_{i_{m}}})^{\frac{1}{i_{m+1}-i_{m}}}$, which is a contradiction.

Case (2): If $l<m$, then $x>(\frac{\delta_{i_{m}}}{\delta_{i_{l}}})^{\frac{1}{i_{m}-i_{l}}}$. By Lemma \ref{lem1} we have
\[
\min\Bigg\{(\delta_{i_{m}})^{\frac{1}{i_{m}}},(\frac{\delta_{i_{m}}}{\delta_{i_{1}}})^{\frac{1}{i_{m}-i_{1}}},\cdots ,(\frac{\delta_{i_{m}}}{\delta_{i_{m-1}}})^{\frac{1}{i_{m}-i_{m-1}}}\Bigg \}=(\frac{\delta_{i_{m}}}{\delta_{i_{m-1}}})^{\frac{1}{i_{m}-i_{m-1}}}.
\]
So $x>(\frac{\delta_{i_{m}}}{\delta_{i_{m-1}}})^{\frac{1}{i_{m}-i_{m-1}}}$, and this is a contradiction. This completes the proof.
\end{proof}

\begin{theorem}\label{th4}
Let $A\in M_n(\mathbb{R}_+)$, where $a_{11}\leq a_{22}\leq \cdots \leq a_{nn}$ and let
\[
 \mathcal{X}_{A}(x)= x^{n}\oplus \delta_{i_{1}}x^{n-i_{1}}\oplus \delta_{i_{2}}x^{n-i_{2}} \oplus \cdots \oplus\delta_{i_{t}}x^{n-i_{t}},
\]
where $1\leq i_1<i_2<\cdots <i_t\leq n$ and where $\delta_{i_j},~1\leq i\leq t$ are  essential terms. 
Then for all $1\leq m\leq t$ we have
\begin{equation}
(\frac{\delta_{i_{m}}}{\delta_{i_{m-1}}})^{\frac{1}{i_{m}-i_{m-1}}}\geq a_{n-i_{m-1},n-i_{m-1}}.
\label{key}
\end{equation}
\end{theorem}

\begin{proof}
Let $1\leq m\leq t$ be given. Now consider the product:
\[ \left(\frac{\delta_{i_{m}}}{\delta_{i_{m}-1}}\right) \left(\frac{\delta_{i_{m}-1}}{\delta_{i_{m}-2}}\right) \cdots \left(\frac{\delta_{i_{m-1}+1}}{\delta_{i_{m-1}}}\right),\]
consisting of $i_m-i_{m-1}$ factors. Using Lemma \ref{l1}, we have $\delta_{i+1}/\delta_i \geq a_{n-i,n-i}$, for $1 \leq i \leq n$. Thus the product above is at least 
\[a_{n-(i_m -1),n-(i_m -1)}a_{n-(i_m -2),n-(i_m -2)} \cdots a_{n-i_{m-1},n-i_{m-1}}.\] Using the hypothesis on the main diagonal entries of $A$  we have that this diagonal product is at least $(a_{n-i_{m-1},n-i_{m-1}})^{i_m-i_{m-1}}$. Thus it follows that 
\begin{eqnarray*}
\left(\frac{\delta_{i_{m}}}{\delta_{i_{m-1}}}\right)&=&\left(\frac{\delta_{i_{m}}}{\delta_{i_{m}-1}}\right) \left(\frac{\delta_{i_{m}-1}}{\delta_{i_{m}-2}}\right) \cdots \left(\frac{\delta_{i_{m-1}+1}}{\delta_{i_{m-1}}}\right), \\
&\geq&  (a_{n-i_{m-1},n-i_{m-1}})^{i_m-i_{m-1}}.
\end{eqnarray*}
\end{proof}
\begin{theorem}
\label{trop_cor}
Let $A\in M_{n}(\mathbb{R}_{+})$ and let $1\leq k<n$. Then  $\displaystyle conv_{\otimes}\left(\sigma_{trop}^{k}(A)\right)\subseteq W_{max}^{k}(A)$.
\end{theorem}

\begin{proof}
 Using (\ref{p_eq}), we may assume without loss of generality that $a_{11}\leq\cdots\leq a_{nn}$.   Following the discussion before Remark \ref{r19}, it follows  that the characteristic maxpolynomial of $A$ can be written as
 \begin{eqnarray*}
 \mathcal{X}_{A}(x)&=& x^{n}\oplus \delta_{i_{1}}x^{n-i_{1}}\oplus \delta_{i_{2}}x^{n-i_{2}} \oplus \cdots \oplus\delta_{i_{t}}x^{n-i_{t}}\\
 &=& (x\oplus (\delta_{i_1})^{\frac{1}{i_1}})^{i_1}(x\oplus (\frac{\delta_{i_2}}{\delta_{i_1}})^{\frac{1}{i_2-i_1}})^{i_2-i_1}\ldots (x\oplus (\frac{\delta_{i_t}}{\delta_{i_{t-1}}})^{\frac{1}{i_t-i_{t-1}}})^{i_t-t_{t-1}}x^{n-i_{t}},
 \end{eqnarray*}
 where $i_{0}=1\leq i_{1}<i_{2}<\cdots <i_{t}\leq i_{t+1}=n$  and where $\delta_{i_j},~1\leq j\leq t$ are the essential terms. By Lemma \ref{orderroots}, 
 $
 (\delta_{i_1})^{\frac{1}{i_1}}\geq  (\frac{\delta_{i_2}}{\delta_{i_1}})^{\frac{1}{i_2-i_1}}\geq \ldots  \geq (\frac{\delta_{i_t}} {\delta_{i_{t-1}}})^{\frac{1}{i_{t}-i_{t-1}}}
$
  are the nonzero algebraic max eigenvalues of $A$ with corresponding multiplicities $i_1,i_2-i_1,\ldots ,i_t-i_{t-1}$, respectively.  So we have
  \[
  0\leq n-i_{t}<n-i_{t-1}<\cdots <n-i_{1}\leq n-1.
  \]
If $n-i_{s}<k\leq n-i_{s-1},~ 1\leq s\leq t+1,$ then we have
 \[
 \sigma_{trop}^{k}(A)=\displaystyle\left\lbrace\displaystyle\left(\delta_{i_{1}}\right)^{\frac{1}{i_{1}}},   \displaystyle\left(\frac{\delta_{i_{2}}}{\delta_{i_{1}}}\right)^{\frac{1}{i_{2}-i_{1}}}, \ldots,  \displaystyle\left(\frac{\delta_{i_{s}}}{\delta_{i_{s-1}}}\right)^{\frac{1}{i_{s}-i_{s-1}}}\right\rbrace.
 \]
 By Theorem \ref{th4}
 \[
 \displaystyle\min_{z\in \sigma_{trop}^{k}(A)}z=\displaystyle\left(\frac{\delta_{i_{s}}}{\delta_{i_{s-1}}}\right)^{\frac{1}{i_{s}-i_{s-1}}}\geq a_{n-i_{s-1}, n-i_{s-1}}\geq a_{kk}.
 \]
 Since
 \[
 \displaystyle\max_{z\in \sigma_{trop}^{k}(A)}z= \displaystyle\max_{z\in \sigma_{trop}(A)}z\leq d,
 \]
where $d=\displaystyle\max_{1\leq i,j\leq n} a_{ij}$ 
 and since $W_{max}^{k}(A)=[a_{kk}, d]$ by Theorem \ref{Main Theorem},  the proof is complete.
 
\end{proof}

%
%



\section{\textbf{Corrections of \cite[Section 5]{TZAH}}}

Let $ A \in M_{n}(\mathbb{R}_{+})$ and $c=[c_{1}, c_{2}, \ldots, c_{n}]^{t} \in \mathbb{R}_{+}^{n}.$  
In \cite[Section 5]{TZAH} we defined the {\em max $c-$numerical range of $A$} as follows:
\[
W_{max}^{c}(A)=\{\bigoplus_{i=1}^{n} c_{i}(x^{(i)})^{t}\otimes A \otimes x^{(i)}:~X=[x^{(1)}, x^{(2)}, \ldots, x^{(n)}] \in M_{n}(\mathbb{R}_{+}), ~X \in \mathcal{U}_{n}  \}.
\]
Clearly $W_{max}^{c}(A)=\{tr_{\otimes}(C\otimes X^{t}\otimes A \otimes X): X \in \mathcal{U}_{n} \},$  where
$C=diag(c_{1}, \ldots, c_{n}),  c=[c_{1}, c_{2}, \ldots, c_{n}]^{t} \in \mathbb{R}_{+}^{n}.$  

 It was stated in \cite[Section 5]{TZAH} that
\begin{equation}
W_{max}^{c}(A)=\{ c_{k}(\oplus_{i=1}^{n}a_{ii}): k=1,2, \ldots, n \},
\label{prva}
\end{equation}
and 
\begin{equation}
\label{druga}
conv_{\otimes}(W_{max}^{c}(A))=[\min_{1 \leq k \leq n}c_{k}(\oplus_{i=1}^{n}a_{ii}),~ \oplus_{k=1}^{n}c_{k}(\oplus_{i=1}^{n}a_{ii})].
	\end{equation}
However, equations (\ref{prva}) in (\ref{druga}) are not true in general. 	
In fact, we have 
\begin{equation}
W_{\max}^c(A)=\{\displaystyle\oplus_{i=1}^nc_ia_{\sigma(i),\sigma (i)}:\sigma\in \sigma_n\}.
\label{tretja}
\end{equation}
Indeed,  suppose $\lambda\in W_{\max}^c(A)$.  Hence there exists $X=[x^{(1)},\dots ,x^{(n)}]\in\mathcal{U}_n$ such that $\lambda=\displaystyle\oplus_{i=1}^n c_i\otimes (x^{(i)})^t\otimes A\otimes x^{(i)}$. Since $X\in\mathcal{U}_n$   there exists  $\sigma\in \sigma_n$ such that 
\[
x^{(i)}=e_{\sigma (i)}, 1\leq i\leq n,
\]
where $e_{\sigma (i)}=[0,\dots ,0,1,0,\dots ,0]^t\in\mathbb{R}_+^n.$
Therefore 
\[
\lambda =\oplus_{i=1}^nc_ie_{\sigma (i)}^t\otimes A\otimes e_{\sigma (i)}=\oplus_{i=1}^nc_ia_{\sigma (i),\sigma (i)},
\]
which establishes (\ref{tretja}).

 Consequently, if $c_1=\dots =c_n$ or $a_{11}=\dots =a_{nn}$, then $W_{\max}^c(A)$ is a singleton set.  Furtheremore, it follows that 
 \[
conv_{\otimes}(W_{\max}^c(A))=[\displaystyle\min_{\sigma\in S_n}\displaystyle\oplus_{i=1}^nc_ia_{\sigma(i),\sigma (i)},  \displaystyle\max_{\sigma\in S_n}\displaystyle\oplus_{i=1}^nc_ia_{\sigma(i),\sigma (i)}]. 
 \]
 
 
 \bigskip
 It was correctly noted in  \cite[Example 5]{TZAH} that $conv_{\otimes}(W_{max}^{C}(A))=[\displaystyle\min_{1\leq i\leq n}a_{ii}, \oplus_{i=1}^{n}a_{ii}],$ where
 $A=(a_{ij})\in M_{n}(\mathbb{R}_{+})$ and
\[
C=\left [
\begin{array}{cccc}
1&0&\cdots&0\\
0&0&\cdots&0\\
\vdots &\vdots &\cdots&\vdots\\
0&0&\cdots&0
\end{array}
\right]\in M_{n}(\mathbb{R}_{+}).
\] 
We extend this idea in the following remark.
\begin{remark}\label{extend_C numerical range}
Let $C=(c_{ij})\in M_{n}(\mathbb{R}_{+})$ such that $c_{rs}=1$ and $c_{ij}=0$ elsewhere. Since
 \begin{equation*}
W_{max}^{C}(A)=  
  \begin{cases}
   \{a_{ii}: 1\leq i\leq n\}  &  r=s  \\ 
    \{a_{ij}: i, j\in\{1, 2, \ldots, n\},~ i\neq j\} &  r\neq s,
\end{cases}
\end{equation*}
one has
\begin{equation*}
conv_{\otimes}(W_{max}^{C}(A))=  
  \begin{cases}
   [\displaystyle\min_{1\leq i\leq n}a_{ii}, \oplus_{i=1}^{n}a_{ii}]  &  r=s  \\ 
    [\displaystyle\min_{1\leq i, j\leq n, i\neq j}a_{ij}, \oplus_{i, j=1,~ i\neq j}^{n}a_{ij}] &  r\neq s.
\end{cases}
\end{equation*}
 
\end{remark}
Most of the following result was correctly stated in \cite[Theorem 5, properties (i)-(iv), (vi)]{TZAH}. There was a typing error in  \cite[Theorem 5, property (v)]{TZAH}, which we correct below.
We include details of the proof.
\begin{theorem}\label{THCNumerical}
Let $A, C \in M_{n}(\mathbb{R}_{+})$. Then the following assertions hold:
\begin{itemize}
\item[(i)]  $W^{C}_{max}(\alpha A\oplus  \beta I_{n}) = \alpha W^{C}_{max}(A) \oplus  \beta tr_{\otimes}(C),$ where $\alpha, \beta \in  \mathbb{R}_{+};$

\item[(ii)] $ W_{\max }^{C}(A\oplus B) \subseteq W_{\max }^{C}(A)\oplus W_{\max }^{C}(B)$ and $W_{\max }^{C\oplus D}(A) \subseteq W_{\max }^{C}(A)\oplus W_{\max }^{D}(A),$ where $B,D \in M_{n}(\mathbb{R}_{+});$

\item[(iii)]   $W_{\max }^{C}(U^{t}\otimes A\otimes U)= W_{\max }^{C}(A),$  where  $U\in U_{n};$

\item[(iv)]  If  $C^{t}=C,$ then    $W_{\max }^{C}( A^{t})= W_{\max }^{C}(A);$ 

\item[(v)]    If  $C=\alpha I_{n},$  where $\alpha \in \mathbb{R}_{+},$  then  $W_{\max }^{C}(A)=\{ \alpha tr_{\otimes}(A) \}$ 

\item[(vi)]    $W_{\max }^{ C}( A)=W_{\max }^{ A}( C)$.
\end{itemize}
\end{theorem}
\begin{proof}
\begin{itemize}
\item[(i)] Let $z\in W_{\max }^{ C}(\alpha A\oplus \beta I_{n})$. Then $z=\displaystyle  tr_{\otimes}\left(C\otimes X^{t}\otimes (\alpha A\oplus \beta I_{n})\otimes X \right)$  for some $X\in \mathcal{U}_{n}$ and hence $z=\alpha \displaystyle tr_{\otimes} (C\otimes X^{t} \otimes A \otimes X)\oplus \beta tr_{\otimes}(C)$. This implies that $z\in \alpha W_{max}^{C}(A)\oplus \beta tr_{\otimes}(C)$. For the reverse inclusion, let $z\in \alpha W_{max}^{C}(A)\oplus \beta tr_{\otimes}(C)$. So $z=\alpha \displaystyle\left(tr_{\otimes} (C\otimes X^{t} \otimes A \otimes X)\right)\oplus \beta tr_{\otimes}(C)$ for some $X\in \mathcal{U}_{n}$ and it follows that $z\in W_{max}^{C}(\alpha A\oplus \beta I_{n})$.

\item[(ii)]  Let $z\in W_{\max }^{ C}(A\oplus B)$. Then $z=tr_{\otimes}\displaystyle \left(C\otimes X^{t} \otimes (A\oplus B) \otimes X\right)$ for some $X\in \mathcal{U}_{n}$ and hence 
$z=tr_{\otimes}(C\otimes X^{t} \otimes A\otimes X)\oplus tr_{\otimes}(C\otimes X^{t} \otimes B\otimes X)$ for some $X\in \mathcal{U}_{n}$. This implies that 
$z\in W_{\max }^{C}(A)\oplus W_{\max }^{C}(B)$. For a proof of second part, let $z\in W_{\max }^{C\oplus D}(A) $. So $z=tr_{\otimes}\displaystyle \left((C\oplus D)\otimes X^{t} \otimes A\otimes X\right)$ for some $X\in \mathcal{U}_{n}$ and hence 
$z=tr_{\otimes}\displaystyle (C\otimes X^{t} \otimes A\otimes X)\oplus tr_{\otimes}\displaystyle (D\otimes X^{t} \otimes A\otimes X)$. This implies that $z\in W_{\max }^{C}(A)\oplus W_{\max }^{D}(A)$.

\item[(iii)] Let $z\in W_{\max }^{C}(U^{t}\otimes A\otimes U)$. Then $z=\displaystyle  tr_{\otimes}\left(C\otimes X^{t}\otimes (U^{t}\otimes A\otimes U)\otimes X \right)$ for some $X\in \mathcal{U}_{n}$ and hence $z=\displaystyle  tr_{\otimes}\left(C\otimes (U\otimes X)^{t}\otimes A\otimes (U\otimes X) \right)$ for some $X\in \mathcal{U}_{n}$. Since $U\otimes X \in \mathcal{U}_{n},$ one has $z\in W_{\max }^{C}(A)$. For the reverse inclusion, let $z\in W_{\max }^{C}(A)$. Thus $z=tr_{\otimes}\displaystyle \left(C\otimes X^{t} \otimes A\otimes X\right)$ for some $X\in \mathcal{U}_{n}$. Set $U^{t}\otimes A\otimes U=B,$ or $A=U\otimes B\otimes U^{t}$. Therefore
\begin{eqnarray*}
z &=&  tr_{\otimes}\displaystyle \left(C\otimes X^{t} \otimes U\otimes B\otimes U^{t}\otimes X\right) \\
&=& tr_{\otimes}\displaystyle \left(C\otimes Y^{t}\otimes B\otimes Y\right)\\
&\in & W_{\max }^{C}(B)=W_{\max }^{C}(U^{t}\otimes A\otimes U),
\end{eqnarray*}
where $Y=U^{t}\otimes X\in \mathcal{U}_{n}$.
\item[(iv)]  If $C^{t}=C,$  then 
\begin{eqnarray*}
W_{\max }^{C}( A^{t})&=&\{tr_{\otimes}(C\otimes X^{t}\otimes A^{t}\otimes X): X\in \mathcal{U}_{n}\}\\
&=&\{tr_{\otimes}(C^{t}\otimes X^{t}\otimes A^{t}\otimes X): X\in \mathcal{U}_{n}\}\\
&=& \{tr_{\otimes} (X^{t}\otimes A^{t}\otimes X\otimes C^{t}): X\in \mathcal{U}_{n}\}\\
&=& \{tr_{\otimes} (C\otimes X^{t}\otimes A\otimes X): X\in \mathcal{U}_{n}\}\\
&=& W_{\max }^{C}( A).  
\end{eqnarray*}

  \item[(v)]  If $C=\alpha I_n$, then 
  \begin{eqnarray*}
W_{\max }^{C}( A)&=&\{tr_{\otimes}(C\otimes X^{t}\otimes A\otimes X): X\in \mathcal{U}_{n}\}\\
&=&\{tr_{\otimes}( \alpha I_{n}\otimes X^{t}\otimes A\otimes X): X\in \mathcal{U}_{n}\}\\
&=& \{\alpha tr_{\otimes} ( X^{t}\otimes A\otimes X): X\in \mathcal{U}_{n}\}\\
&=& \{\alpha tr_{\otimes} (A)\}.
\end{eqnarray*}

  \item[(vi)] Finally, 
\begin{eqnarray*}
W_{\max }^{C}( A)&=&\{tr_{\otimes}(C\otimes X^{t}\otimes A\otimes X): X\in \mathcal{U}_{n}\}\\
&=&\{tr_{\otimes}(A\otimes X\otimes C\otimes X^{t}): X\in \mathcal{U}_{n}\}\\
&=& W_{\max }^{ A}( C).  
\end{eqnarray*}
\end{itemize}
\end{proof}

\bigskip

{\bf Acknowledgments.} 

The authors thank S. Gaubert for pointing out and correcting some mistakes from \cite{TZAH}.

The third author acknowledges a partial support of the Slovenian Agency for Research and Innovation (grants P1-0222, J1-8133, J2-2512 and J1-8155). 

Fallat's research is supported in part by an NSERC Discovery Grant, RGPIN--2019--03934.

\bibliographystyle{amsplain}

\end{document}